%% file: Bayes_JS_SIAMUQ.tex
\documentclass[review,onefignum,onetabnum]{siamonline190516}

\input{ex_shared2}

\ifpdf
\hypersetup{
  pdftitle={{Empirical Bayesian Inference Using Joint Sparsity}}
  pdfauthor={J. Zhang, A. Gelb, and T. Scarnati}
}
\fi

\graphicspath{{./figures/}{./figures/individual_figs/}}

\begin{document}

\maketitle

\begin{abstract}
\input{abstract}
\end{abstract}


\begin{AMS}
 {62C12,65C40,68U10}
\end{AMS}

\section{Introduction}
\label{sec:Intro}
\input{introduction}

\section{The Bayesian approach for solving inverse problems}
\label{sec:bayes}
\input{bayes}

\section{Evaluating the Posterior}
\label{sec:posterior}
\input{posterior}

\section{Numerical Experiments}
\label{sec:numerics}

\input{numerics}

\section{Conclusions}
\label{sec:conclusions}
\input{conclusions}

\vspace{5mm}
\noindent {\bf Acknowledgments.}
Anne Gelb's work is supported in part by the NSF grants  DMS \#1502640 and DMS \#1912685,  AFOSR grant \#FA9550-18-1-0316, and ONR MURI grant \#N00014-20-1-2595.
Theresa Scarnati's work is supported in part by AFOSR LRIR \#18RYCOR011.

\appendix
\section{Polynomial Annihilation (PA) Matrix Construction}
\label{appendix:PA}
\input{appendixPA}


\bibliographystyle{siamplain}
\bibliography{references}
\end{document}

%% file: ex_shared2.tex
\usepackage{lipsum}
\usepackage{amsfonts}
\usepackage{graphicx}
\usepackage{epstopdf}
\usepackage{algorithmic}
\usepackage[margin=1in]{geometry}
\usepackage{amssymb}
\usepackage{amsmath}
\usepackage{graphicx}
\usepackage{url}
\usepackage{color}
\usepackage{framed}
\usepackage{algorithm}
\usepackage{mathrsfs}
\usepackage{bm}
\usepackage[normalem]{ulem}
\usepackage{multirow}
\ifpdf
  \DeclareGraphicsExtensions{.eps,.pdf,.png,.jpg}
\else
  \DeclareGraphicsExtensions{.eps}
\fi

\usepackage{enumitem}
\setlist[enumerate]{leftmargin=.5in} 
\setlist[itemize]{leftmargin=.5in}

\newsiamremark{remark}{Remark}
\newsiamremark{hypothesis}{Hypothesis}
\crefname{hypothesis}{Hypothesis}{Hypotheses}
\newsiamthm{claim}{Claim}
\newsiamthm{example}{Example}

\title{Empirical Bayesian Inference using Joint Sparsity}

\author{Jiahui Zhang\thanks{Department of Mathematics, University of Minnesota, Minneapolis, MN 
  (\email{zhan.7554@umn.edu}
  ).}
\and Anne Gelb\thanks{Department of Mathematics, Dartmouth College, Hanover, NH
  (\email{annegelb@math.dartmouth.edu}).}
\and Theresa Scarnati\thanks{Qualis Corporation, Huntsville, AL
  (\email{tscarnati@qualis-corp.com}).}}

\usepackage{hyperref}
\usepackage{amsopn}

\newcommand{\AG}{\textcolor{red}}

\newcommand*{\bbR}{\mathbb{R}}
\newcommand*{\supp}{\mathrm{supp}}
\newcommand*{\argmin}[1]{{\underset{#1}{\operatorname{argmin}}}}
\newcommand*{\argmax}[1]{{\underset{#1}{\operatorname{argmax}}}}

\newcommand{\bes}[1]{\begin{equation*} #1 \end{equation*}}
\newcommand*{\nm}[1]{{\left\|#1\right \|}}

%% file: abstract.tex
This paper develops a new empirical Bayesian inference algorithm for solving a linear inverse problem given multiple measurement vectors (MMV) of under-sampled and noisy observable data.  Specifically, by exploiting the {\em joint sparsity} across the multiple measurements in the sparse domain of the underlying signal or image,  we construct a new {{\em support informed}} sparsity promoting prior. Several applications can be modeled using this framework, and as a prototypical example we consider reconstructing an image from synthetic aperture radar (SAR) observations using nearby azimuth angles. 
Our numerical experiments demonstrate that using this new prior not only improves accuracy of the recovery, but also reduces the uncertainty in the posterior when compared to standard sparsity producing priors. 


%% file: introduction.tex
We consider the linear inverse problem modeled as {$\mathbf{Y} = A \mathbf{X} + \mathbf{E}$}  {where $\mathbf{X}, \mathbf{Y}, \mathbf{E}$ are random variables defined over the common probability space $(\Omega, \mathcal{F}, \mathbb{P})$,} {and $\mathbf{X}$ and $\mathbf{E}$ are assumed to be independent.}
Here $\mathbf{X}$ represents the underlying unknown we seek to recover, $\mathbf{Y}$ is the observable data, {$A \in \mathbb{C}^{n \times n}$} is a {deterministic} and known forward operator, and {$\mathbf{E}\sim\mathcal{CN}(0,\sigma^2\mathbb{I}_{{n}})$ models the additive complex Gaussian noise.}   The realizations of the observable data $y_j = Ax + e_j$, $j = 1,\dots,J$, are obtained as a collection of {$J > 1$}  measurements,  and are often referred  to as multiple measurement vectors (MMVs), \cite{MMVChenHuo}. Many real world phenomena may be modeled in this way, {and our study is motivated by applications in synthetic aperture radar (SAR) imaging.  In this case, the observable data are multiple sets of phase history data,  $A$ is a (non-uniform) Fourier transform matrix, and we seek to recover a complex  SAR image  ${\mathbf X}$.}  To simplify our presentation and analysis, our investigation is limited to one-dimensional real-valued signals. It is straightforward to extend our approach to multi-dimensions as well as complex-valued signals and images, however. Additionally, our approach is conducive to parallelizable algorithms. 

A common assumption made in signal\footnote{We use the terms signal, image, and function interchangeably throughout our manuscript.} recovery, which we also use in this investigation, is that the underlying signal is sparse in some domain (e.g.~the gradient or edge domain). In such cases {\em compressive sensing} (CS) algorithms, \cite{Candes, Donoho}, are frequently employed, as they are specifically designed to exploit these sparsity assumptions.  Constructed as a convex optimization problem, a CS algorithm seeks to minimize the $\ell_2$ norm of the data fidelity regularized by the $\ell_1$ norm of the solution in its sparse domain. {For example, given a particular observation $y_j = Ax + e_j$, the solution is obtained as
\begin{equation}\label{eq:l1reg}
x^\ast
   = \argmin{x \in \mathbb{R}^n} \Big\{ \lambda||\mathcal{L}x||_1 + \frac{1}{2} ||y_j-Ax||_2^2 \Big\},
\end{equation}
where $\mathcal{L}$ is a pre-determined sparsifying transform and the regularization parameter $\lambda>0$ is chosen to balance each term's contribution.  For presumably less error in the fidelity term, $||y_j-Ax||_2^2$, a smaller $\lambda$ is chosen and vice versa.}
In this context the $\ell_1$ norm is typically viewed as a surrogate for the $\ell_0$ pseudo-norm since using the latter yields an intractable problem, {\cite{CompressiveSensingCandes}}.  The additional information from the multiple observations collected in the MMVs can be further exploited by employing the {\em joint sparsity} assumption. In a nutshell, joint sparsity leverages the assumption that the sparse domain of the underlying signal should be similar across all collected measurements for that same signal.

Many techniques have been introduced  to improve the robustness and efficacy of CS algorithms in both the single and multiple measurement cases, \cite{candes2008enhancing, daubechies2010iteratively, VBJSGelb, VBJS, VBJS2}.  However, since the recovery is limited to point estimate solutions, it is not possible to quantify the uncertainty. This disadvantage is non-trivial because in practice it is often crucial to know how reliable the recovery is, especially in the case of noisy or limited data availability within in each measurement vector. Moreover, even within the same application, CS algorithms inevitably require hand tuning {of the regularization parameter}, with small changes often resulting in drastically different recovery outcomes.

{Due to these inherent limitations, the Bayesian approach, \cite{kaipio2006statistical}, is often used to recover the posterior distribution of the unknown random variable ${\bf X}$ with respect to ${\mathbf Y}$ as related by the model ${\mathbf Y} = A {\mathbf X} + {\mathbf E}$.}   We note that from the Bayesian perspective, the CS methodology may be viewed as obtaining a maximum a posteriori  (MAP) estimate \cite{murphy2012machine}, which represents a  mode of the posterior distribution. However, the MAP estimate of the posterior distribution may not always provide the {most descriptive approximation of the true underlying signal and furthermore does not yield the required} information about the posterior distribution to quantify the uncertainty of the recovery. Therefore this investigation aims to recover the posterior distribution representing the unknown $\mathbf{X}$ using sampling methods. 

In the Bayes formulation, the posterior density of the unknown $\mathbf{X}$ is proportional to the product of the likelihood and prior densities, with the likelihood corresponding to the fidelity term in a convex optimization problem and the prior to the regularization term. Determining an appropriate prior to quantify assumptions about the unknown signal {is well known to be difficult}. To this end, the sparsity assumption suggests that the $\ell_1$ (Laplace) prior {may be} suitable for this purpose, {and indeed the corresponding MAP estimate results in the same convex CS-based optimization problem in (\ref{eq:l1reg}).}  Although effective in many applications, the standard Laplace prior does not really include all of the prior beliefs about the underlying image.  Specifically, no consideration is given to the {\em locations} of the support within the sparse domain.  We note that this is also the case in standard CS algorithms, as the regularization parameter is constant throughout the sparsity domain.  

This issue has been addressed in CS algorithms using various forms of {\em weighted} $\ell_1$ regularization, \cite{daubechies2010iteratively, candes2008enhancing, VBJSGelb,VBJS}.  In all cases, the {\em spatially varying} weighting parameters are designed to promote solutions that are sparse only in regions that are {likely to be}  sparse, as opposed to promoting {\em global} sparsity.  Unfortunately, many of the proposed re-weighted $\ell_1$ (and $\ell_2$) regularization algorithms are iterative in nature and  sensitive to noise. A review of the advantages and shortcomings of some of these algorithms may be found in \cite{churchill2018edge}. 

In \cite{VBJSGelb,VBJS}, where the MMV case was specifically considered, the  spatially varying weighting parameter was designed by exploiting the {\em joint sparsity} across the measurements in the sparsity domain.  {In particular, the regions of support in the sparse domain were determined to be those in which there were large variances in the measurements.  By contrast, small variances in the measurements indicated sparse regions. The spatially varying weights for the weighted $\ell_1$ (or $\ell_2$) regularization were then constructed accordingly.} The technique, coined {\em the variance based joint sparsity} (VBJS) method, proved to be efficient and robust in recovering images from noisy and under-sampled data. {The method can also be used for a single set of measurements by processing it multiple times to obtain MMVs, \cite{VBJS2}.}

Inspired by these results, in the empirical Bayesian approach used in this investigation,  we introduce a new prior that is constructed by exploiting the joint sparsity assumption.   As is typically done, our prior assumes that the signal is sparse in the sparse domain. But {in our new {\em support informed} sparse prior, we also incorporate the likely regions of support in the sparse domain, thereby  more accurately describing our prior beliefs about the underlying signal}. As in the VBJS method, our technique uses the variance of the measurements in the sparsity domain to determine regions of support. {We then use this information to construct a binary edge mask for the prior.  That is, the weighting parameter in the prior} is binary with a zero value assigned whenever an edge is detected and one otherwise.\footnote{Binary edge masks were also used in the edge-adaptive $\ell_2$ regularization image reconstruction method  developed in \cite{churchill2018edge}.}  In this way we obtain a truly hierarchical prior, which will be described in more detail throughout Section \ref{sec:bayes}.  Since our new method uses observable data to inform the prior it may be categorized  as an empirical Bayesian inference technique. 


For context, we can compare our method to another empirical Bayesian inference technique, {\em sparse Bayesian learning} (SBL), \cite{tipping2001sparse, wipf2004sparse}.  In SBL, one seeks to \textit{learn} the distributions of the hyperparameters of the prior.  For a conjugate prior (given in \cite{tipping2001sparse,wipf2004sparse} by a Gaussian prior with Gamma distributed hyperparameters), these hyperparameters can be approximated using expectation maximization (EM) \cite{dempster1977maximum}. SBL with EM has the same effect as defining a sparse prior that is additionally informed by the supports in the sparse domain. In this regard it is important to point out that (i) SBL is designed specifically to recover sparse signals. That is, the sparse prior is in the domain of recovery.  By contrast our new method is more general, since it also considers priors that are sparse in a transformed space;  (ii)  EM is not computationally efficient for higher dimensional data. While other more efficient methods may be used to learn the hyperparameters, the resulting accuracy is considerably reduced, \cite{wipf2004sparse}; {and (iii)} due to the use of conjugate priors, the density yields a closed form representation of the posterior distribution, and thus sampling can be done using the Gibbs sampling method, {an algorithm} in the family of Markov Chain Monte-Carlo (MCMC) techniques. {Our proposed technique, on the other hand, is not restricted to using only conjugate priors.}  Consequently, the posterior distribtion does not have a closed form solution, which means that we cannot use the Gibbs sampling method. Hence instead we employ the celebrated Metropolis Hastings (MH) MCMC algorithm, \cite{MCMC,MCMCR}. For one dimensional problems, the efficiency of the MH algorithm is comparable to Gibbs sampling. In higher dimensions, Gibbs sampling is inevitably a faster sampling technique. The evaluation of the posteriors derived in this paper is easily parallelized across each dimension, however, so that for many applications enforcing conjugate priors to enable Gibbs sampling does not seem necessary. 

Our new method may also be compared those of \cite{BardsleyUQ,calvetti2020hybrid, calvetti2020sparse}, which in sequel we will refer to  as the {hierarchical Gaussian prior (HGP) approach. In contrast to SBL, HGP techniques are not limited to the recovery of sparse signals, as long as the sparsifying transform in the prior is symmetric positive definite.  This will ensure a closed form solution for the posterior which in turn allows Gibbs sampling to be used.}\footnote{{Based on its derivation, it seems that EM could also be used in the case of SPD sparsifying transforms. This would make SBL applicable to a broader range of signals (not only sparse signals).  While such applications should be considered in future investigations, it does not negate SBL's issue with implementation cost in multi-dimensions.}}  {As in SBL, the prior in HGP  is typically Gaussian with} {Gamma} {distributed  hyperparameters. Unlike SBL, these hyperparameters are not approximated by EM.  In order to maintain a closed form solution for the posterior, the hyperparameters are {\em scalar}, that is, the prior is not informed by the signal's support in the sparsity domain.  In this regard, a point estimate of the HGP is more similar to the original CS approach which does not use spatially varying weighted $\ell_1$ regularization.  Of course HGP has the benefit of being able to quantify the uncertainty of the solution.} We note again that in our new approach we use a hierarchical support informed prior, where the support in the sparse domain is determined {by exploiting the joint sparsity of the multiple measurements in the sparse domain.}



{The rest of this paper is organized as follows.  In Section \ref{sec:bayes} we review the Bayesian approach for solving inverse problems and introduce our new support informed sparse prior.    We discuss how the resulting posterior is evaluated via the MH MCMC algorithm in Section \ref{sec:posterior}.  Numerical experiments used to evaluate our method are performed in Section \ref{sec:numerics}.  Finally, Section \ref{sec:conclusions} provides some concluding remarks.}

%% file: bayes.tex
We consider the linear forward model of the form
\begin{equation}
\label{eq:ForwardModel}
    \mathbf{Y} = A \mathbf{X}+ \mathbf{E},
\end{equation}
{where $\mathbf{X}, \mathbf{Y}, \mathbf{E}$ are random variables 
 defined {over} a common probability space $(\Omega, \mathcal{F}, \mathbb{P})$,} {and $\mathbf{X}$ and $\mathbf{E}$ are assumed to be mutually independent.}
Here $\mathbf{X}$ represents the underlying unknown we seek to recover, $\mathbf{Y}$ is the observable data, {$A \in \mathbb{C}^{n \times n}$} is a known forward operator, and {$\mathbf{E}\sim\mathcal{CN}(0,\sigma^2\mathbb{I}_{{n}})$ models additive complex Gaussian noise.} 
A particular observation {$y_j$ of (\ref{eq:ForwardModel}), $j = 1,\dots J$, where $J$ is the given number of MMVs,}  is given by 
\begin{equation}
\label{eq:PointEstimate}
y = Ax + e,
\end{equation}
{where we have dropped the subscript $j$ to avoid cumbersome notation.}
As noted in the introduction, sparsity based (CS) reconstruction algorithms are often used to obtain the optimal solution to \eqref{eq:PointEstimate} as a point estimate solution of \eqref{eq:ForwardModel}.  {Also as discussed in the introduction, exploring the entire posterior probability of $\bm{X}$ is valuable in assessing the uncertainty of the our recovery.} This can be achieved within a Bayesian framework, which we briefly review below. A complete overview of Bayesian inverse problems may be found in \cite{kaipio2006statistical}.

\subsection{{Bayes' Theorem}}\label{sec:bayestheorem}

Let $x,y,e \in {\mathbb{R}^n}$ and $\lambda\in[0,1]$ be realizations of the random variables $\mathbf{X}$, $\mathbf{Y}$, $\mathbf{E}$ and $\bm{\lambda}$ respectively. Bayes' theorem is given by
\begin{equation}
\label{eq:full_posterior}
   f_{\mathbf{X},\bm{\lambda}|\mathbf{Y}}(x,\lambda|y) =
    \frac{f_{\mathbf{Y|X}}(y|x){f}_{\mathbf{X}| \bm{\lambda}}(x| \lambda)f_{\bm{\lambda}}(\lambda)}{f_{\mathbf{Y}}(y)},
\end{equation}
where  ${f}_{\mathbf{X},\bm{\lambda}|\mathbf{Y}}(x,\lambda|y)$ is the {posterior density function} we seek to recover,  $f_{\mathbf{Y}|\mathbf{X}}(y|x)$ is the {likelihood function}, and $\tilde{f}_{\mathbf{X},\bm{\lambda}}(x,\lambda) = {f}_{\mathbf{X}| \bm{\lambda}}(x| \lambda)f_{\bm{\lambda}}(\lambda)$ is the {prior probability density function}. The prior is broken into two parts, the conditional prior ${f}_{\mathbf{X}| \bm{\lambda}}(x| \lambda)$ and the hyper-prior $f_{\bm{\lambda}}(\lambda)$. Further, $f_{\mathbf{Y}}(y)$, often called the {evidence}, is given by the total law of probability as  
\begin{equation}\label{eq:MarginalDensity}
    f_{\mathbf{Y}}(y) = 
    \int_{\mathbb{R}^n}\int_{[0, 1]}f_{\mathbf{X}, \mathbf{Y}, \bm{\lambda}}(x, y, \lambda) d\lambda dx= 
    \int_{\mathbb{R}^n}\int_{[0, 1]} f_{\mathbf{Y|X}}(y|x){f}_{\mathbf{X}| \bm{\lambda}}(x| \lambda)f_{\bm{\lambda}}(\lambda)d\lambda dx.
\end{equation}

It is reasonable to assume that \eqref{eq:MarginalDensity} is nonzero because {otherwise} the observation $y \in \mathbb{R}^n$ has a probability of zero and is therefore irrelevant in {the applications we are considering.  Moreover,} since calculating \eqref{eq:MarginalDensity} is often computationally intractable, particularly in high dimensional domains, a closed form posterior density is often unavailable. Thus sampling methods such as those described in {Section \ref{sec:IntervalEstimation}} are typically used to approximate the un-normalized version of the posterior \eqref{eq:full_posterior} given by 
\begin{equation}
\label{eq:f_hat}
   \hat{f}_{\mathbf{X}, \bm{\lambda}|\mathbf{Y}}(x, \lambda|y) =
    f_{\mathbf{Y|X}}(y|x){f}_{\mathbf{X}| \bm{\lambda}}(x| \lambda)f_{\bm{\lambda}}(\lambda).
\end{equation}
{In our numerical experiments} we use the Metropolis-Hastings (MH) algorithm to draw samples of (\ref{eq:f_hat}). 

In what follows, Section \ref{sec:likelihood} defines the likelihood $f_{\mathbf{Y|X}}(y|x)$ with respect to the proposed linear model (\ref{eq:ForwardModel}), while Section \ref{sec:hyper-prior} describes how we estimate the hyper-prior $\lambda$. Section \ref{sec:prior} focuses on the definition of the conditional prior ${f}_{\mathbf{X}|\bm{\lambda}}(x| \lambda)$. There we will first discuss {how} sparsity enforcing priors are often employed in practice and then introduce our novel support informed sparse prior that is {empirically} constructed by exploiting the joint sparsity across the multiple measurements in the sparse domain of the underlying signal.  Finally, in Section \ref{sec:post} we define and explore the posterior $\hat{f}_{\mathbf{X}, \bm{\lambda}|\mathbf{Y}}(x, \lambda|y)$. 

\subsection{Estimating the likelihood}
\label{sec:likelihood}

Returning to the model in \eqref{eq:ForwardModel}, observe that
\[
f_{\mathbf{Y}|\mathbf{X}, \mathbf{E}}(y|x, e) = \delta(y-Ax-e),
\]
where $\delta(\cdot)$ is the Dirac delta function. Marginalizing with respect to the  noise $e$ yields
\begin{align}
    f_{\mathbf{Y}|\mathbf{X}}(y|x) &= \int_{\mathbb{R}^n} f_{\mathbf{Y}|\mathbf{X}, \mathbf{E}}(y|x, e)f_{\mathbf{E}|\mathbf{X}}(e|x)de\nonumber\\
    &= \int_{\mathbb{R}^n} \delta(y-Ax-e)f_{\mathbf{E}|\mathbf{X}}(e|x)de\nonumber\\
    &= f_{\mathbf{E}|\mathbf{X}}(y-Ax|x)\nonumber\\
    &= f_{\mathbf{E}}(y-Ax),\label{eq:errordensity}
\end{align}
with the last equality resulting from the independence between $\mathbf{X}$ and $\mathbf{E}$, \cite{kaipio2006statistical}. Hence we see that the likelihood function, $f_{\mathbf{Y}|\mathbf{X}}(y|x)$, depends on the noise variable $\mathbf{E}$. Since $\mathbf{E}$ is Gaussian with  mean  $0$ and {covariance} matrix $\Sigma = \sigma^2 \mathbb{I}_{{n}}$, the noise density is given by
\begin{equation}\label{eq:NoiseDensity}
    f_{\mathbf{E}}(e) \propto \exp \{ - \frac{1}{2\sigma^2}(e-0)^T (e-0) \} = \exp \{ - \frac{||e||^2_2}{2 \sigma^2} \}.
\end{equation}
Thus from  \eqref{eq:errordensity} and \eqref{eq:NoiseDensity} we obtain the likelihood density
\begin{equation}\label{eq:LikelihoodDensity}
    f_{\mathbf{Y}|\mathbf{X}}(y|x)=f_{\mathbf{E}}(y-Ax) \propto \exp\{ - \frac{1}{2\sigma^2}||y-Ax||^2_2\}.
\end{equation}

\subsection{Estimating the hyper-prior} \label{sec:hyper-prior}

{We use $K$-fold cross validation to estimate the hyper-prior, \cite{murphy2012machine}.  To summarize this procedure, we begin by approximating}
\begin{align}\label{eq:hyper_prior_posterior1}
    f_{\bm{\lambda}}(\lambda) \approx f_{\bm{\lambda}|\mathbf{Y}}(\lambda|y) = 
    \frac{f_{\mathbf{Y}| \bm{\lambda}}(y | \lambda)\tilde{f}_{\bm \lambda}(\lambda)}{f(y)} \propto f_{\mathbf{Y}|\bm{\lambda}}(y|\lambda)\tilde{f}_{\bm{\lambda}}(\lambda),
\end{align}
where $\tilde{f}_{\bm \lambda}(\lambda)$ is the prior on the hyper-prior. By assuming $\lambda\sim U[0,1]$ so that $\tilde{f}_{\bm \lambda}(\lambda) = 1$, we have
\begin{align}
\label{eq:hyper_prior_posterior2}
    f_{\bm{\lambda}}(\lambda) \propto f_{\mathbf{Y}|\bm{\lambda}}(y|\lambda).
\end{align}
The mode of the distribution (\ref{eq:hyper_prior_posterior2}) can be approximated using the MAP estimate as
\begin{equation}
\label{eq:lambda_hat_MLE}
    \hat{\lambda} 
    = \argmax{\lambda} \{ f_{\mathbf{Y} |\bm{\lambda}}(y|\lambda) \}
    = \argmax{\lambda} \int_{\mathbb{R}^n} f_{\mathbf{Y}|\mathbf{X}}(y|x)f_{\mathbf{X}|\bm{\lambda}}(x|\lambda) dx.
\end{equation}
Integration with respect to $x$ over the space $\mathbb{R}^{n}$ is not computationally feasible. We therefore instead employ $K$-fold cross-validation to approximate $\hat{\lambda}$ given the data $y_j \in \Omega$ for $j = 1,\dots,J$.   The general process for $K$-fold cross-validation is to first calculate $M$ MAP estimates of the unknown from a subset of $M < J$ observable measurements using sample candidates of the hyper-prior $\hat{\lambda}$. These are then used to generate $M$ {\em training vectors}, which are in turn compared to the remaining $J-M$ observable measurements, or {\em testing vectors}. The whole process is then repeated for $K$ independent trials, and 
$\hat{\lambda}$ 
is chosen to minimize the mean square error of (\ref{eq:ForwardModel}) between all training and testing vectors.
The procedure is summarized in Algorithm \ref{alg:Cross-Validation}.
\begin{algorithm}[H]
\caption{$K$-Fold Cross-Validation to Estimate $\hat{\lambda}$.}
\label{alg:Cross-Validation}
\begin{algorithmic}[1]
\STATE{\textbf{Input:} Set of multiple measurements (MMVs) ${\bf y} = [y_1, y_2, y_3, \ldots, y_J]$.}
\STATE{\textbf{Output:} Point Estimate of the hyper-prior $\hat{\lambda}$.}
\FOR{$k = 1$ to $K$}
	\STATE{Randomly partition MMVs into $M$ training vectors and $J-M$ testing vectors.}
	\FOR{$i=1$ to $M$}
		\STATE{Sample candidate hyper-prior from $\tilde{\lambda}_{i,k}\sim U[0,1]$.} 
		\STATE{Calculate the MAP estimate of the unknown as
    			$$ 
    			\hat{x}_{i,k} = \argmax{x} \{f_{\mathbf{X},\bm{\lambda}|\mathbf{Y}}(x,\tilde{\lambda}_{i,k}|y_i)\}.
    $$}
		\STATE{Compute training data according to  $\hat{y}_{i,k} = A\hat{x}_{i,k}.$}
		\STATE{Evaluate the training data by calculating
		 $$ 
		E_{i,k} = \frac{1}{J-M}\sum_{l=M+1}^{J}\text{MSE}(\hat{y}_{i,k},y_l),
		$$
		where $\{y_l\}_{l = M+1}^J$ are the partitioned $J-M$ testing vectors.}
		\ENDFOR
		\ENDFOR
		\STATE{Choose $(i^*,k^*)= \argmin{i,k} \hspace{1mm} E_{i,k}$ and set $\hat{\lambda} = \tilde{\lambda}_{i^*,k^*}$.}
\end{algorithmic}
\end{algorithm}

{The hyper-prior $\hat{\lambda}$ recovered from Algorithm \ref{alg:Cross-Validation} approximates the solution to (\ref{eq:lambda_hat_MLE}) such that  
\begin{equation}
\label{eq:lambda_hat_delta}
    f_{\bm{\lambda}}(\lambda) \approx \delta_{\hat{\lambda}}(\lambda). 
\end{equation}
Intuitively, \eqref{eq:lambda_hat_delta} implies that we are simply concentrating the density (\ref{eq:hyper_prior_posterior2}) on the value $\hat{\lambda}$, which is our estimate of the mode of the hyper-prior distribution (\ref{eq:lambda_hat_MLE}).  {Note that the $K$ fold validation process assumes that the distribution for $\lambda$ is unimodal.    Previous efforts, see e.g. \cite{murphy2012machine}, have also demonstrated that (\ref{eq:lambda_hat_delta}) is a good  approximation of (\ref{eq:lambda_hat_MLE}).}

\begin{remark}
The hyper-prior $\hat{\lambda}$ recovered in Algorithm \ref{alg:Cross-Validation} is consistent with how one ideally chooses the regularization parameter in the CS approximation given by \eqref{eq:l1reg}, or as we will again see in \eqref{eq:MAPEstimate}, where the {fidelity and regularization} parameters are combined to form a single parameter.  Specifically this parameter should be chosen to offset the variance of the noise, and the $K$-fold cross validation method provides an empirical approach to calculating this variance.
\end{remark}

\subsection{A support informed sparsity prior}\label{sec:prior}

The main focus of this investigation is to construct a support informed sparsity prior $\tilde{f}_{\mathbf{X},\bm{\lambda}}(x,\lambda) = {f}_{\mathbf{X}| \bm{\lambda}}(x| \lambda)f_{\bm{\lambda}}(\lambda)$ in \eqref{eq:full_posterior}.  As will be demonstrated in what follows, this will reduce the uncertainty of the posterior for models given by \eqref{eq:ForwardModel}. 

Using the approximation in (\ref{eq:lambda_hat_delta}) we can make the approximation
\begin{equation*}
    \tilde{f}_{\mathbf{X},\bm{\lambda}}(x,\lambda) =  {f}_{\mathbf{X}| \bm{\lambda}}(x| \lambda)f_{\bm{\lambda}}(\lambda) \approx {f}_{\mathbf{X}| \bm{\lambda}}(x| \hat\lambda),
\end{equation*}
and it remains to define the prior through ${f}_{\mathbf{X}| \bm{\lambda}}(x| \hat\lambda)$.
To do so, we make two assumptions: (i) the prior enforces the sparsity of some transformation of $\mathbf{X}$ and (ii) there are multiple observations $y$ of $\mathbf{Y}$ for each realization $x$ of $\mathbf{X}$.\footnote{{Similarly to what was done in \cite{VBJS2, VBJSGelb}, it is possible to relax this second assumption and instead have  a single measurement vector {\em processed} multiple times.  For purposes of brevitity we exclude such considerations in this investigation.}} 

{There are many applications for which it is reasonable to assume that the desired solution is sparse in some domain. In our investigation, the underlying signal is assumed to be piecewise smooth, implying that the edge domain is sparse.  It is important to note that as a departure from the traditional sparse priors commonly used in formulating the posterior, our new support informed sparse prior considers both the {\em existence} of the sparse domain as well as the {\em locations} of the non-zero values within it, i.e., its support.} 

 \subsubsection{{Sparsity promoting priors}} 
 \label{sec:prior_smv}

Within the Bayesian formulation, sparsity {may be} enforced through the Laplace prior given by
\begin{equation}\label{eq:laplace_prior}
    f_{\mathbf{X}| \bm{\lambda}}(x| \hat\lambda) \propto \exp\{ - \hat\lambda ||\mathcal{L}x||_1 \},
\end{equation}
where the operator $\mathcal{L}$ is a transform operator used to approximate the sparse domain.  In this investigation we choose $\mathcal{L}$ to be the {\em polynomial annihilation} (PA) transform, \cite{ArchibaldGelb,VBJS}, {which is designed to approximate the corresponding edge function of a piecewise smooth function when acting on the unknown (given by $x$ in \eqref{eq:laplace_prior}) {at specified grid points} in the spatial domain.\footnote{A short description of the PA transform is provided in Appendix \ref{appendix:PA}.}  We note that our new method is not defined by the choice of sparsifying transform (here ${\mathcal L}$), {nor by the sparsity domain (here the edge domain)}, and that other choices may be more appropriate for different applications.

Regardless of how both the hyper-parameter $\hat{\lambda}$ and sparsifying transform matrix $\mathcal{L}$ are determined, the Laplace prior given in (\ref{eq:laplace_prior}) may not provide the most information about what is known about the underlying signal.  As discussed {already}, a better prior would make use of the support {\em locations} in the sparse domain.  
Therefore, following the ideas in \cite{VBJSGelb,VBJS}, we develop a new support informed sparse prior.  In so doing, we better capture the piecewise smoothness of the underlying signal.  A couple of remarks are in order:

\begin{remark}
{Our numerical experiments confirm that using the PA transform with $m > 1$ is appropriate when constructing the support informed sparse prior for a piecewise smooth function.  In particular using $m = 1$, which is equivalent to total variation (TV), does not perform as well.  In general when the underlying function is {\em not} constant in smooth domains, a {\em higher degree} polynomial annihilation is needed to recover an accurate representation of the sparse domain.   This issue is extensively discussed for the compressive sensing framework in \cite{ArchibaldGelb}, where higher order PA is used to avoid the ``staircasing'' effect caused by TV regularization. Similarly, the high order PA describes a better prior in our approach used here.}
\end{remark}
\begin{remark}
\label{rem:horseshoe}
{To be clear, there are a variety of sparsity promoting priors to choose from, including but not limited to the spike and slab prior and the horseshoe prior \cite{bhadra2019lasso, piironen2017sparsity}.   {Hierarchical Gaussian priors, which provides the basic framework for the methods described in \cite{BardsleyUQ, tipping2001sparse}, may also be used.}   {Regardless of the type of prior chosen, some transform operator ${\mathcal L}$ to the sparse domain is still needed.}  This investigation does not intend to provide a complete survey of the advantages and drawbacks of using various sparsity promoting  priors. Instead we demonstrate that our proposed prior is a departure from the aforementioned priors  because it {\em explicitly} includes the support locations in the sparse domain.  Indeed, such knowledge may be incorporated into some of the aforementioned choices.  Even further, our approach may be used to form mixed priors.  That is, because our spatially adaptive prior encourages the separation of scales in the sparse domain, different priors may be used in different regions of the domain.  The Laplace prior in \eqref{eq:laplace_prior} is mentioned as a prototype, and is {a common} choice when interested in {constructing the MAP estimate of the posterior} (CS approach).  As will be demonstrated by our numerical examples, our support informed sparse prior yields better accuracy and reduced uncertainty than other priors typically used. More information about sparse priors may be found in \cite{piironen2017sparsity}. }
\end{remark}

\subsubsection{Joint sparsity}\label{subsec:jointsparse}

As noted in the introduction, many applications collect multiple measurement vectors (MMVs) of observable data. By exploiting the {\em joint sparsity} in the sparse domain of the underlying signal that can be obtained from the MMVs, a technique was developed in \cite{VBJSGelb, VBJS} to improve the robustness of the MAP estimate.\footnote{{We note that the method was derived from the CS perspective.}} For reasons that will become apparent later in this section, the technique was coined the Variance Based Joint Sparsity (VBJS) method. We now adapt VBJS  to construct a new {\em support informed sparse prior}. 
We then demonstrate that this new prior yields reduced uncertainty and improved efficiency in the {sampling of}  the  posterior in \eqref{eq:f_hat} when compared to other typical prior definitions, such as the Laplace prior in (\ref{eq:laplace_prior}).  

As previously mentioned, the domain of the problem is $\mathbb{R}^n$ so that the unknown $x$ may be represented by a $n$ length vector. {A $d$-dimensional problem will then yield an $n^d$ length vector representation of the unknown, and as already noted our method is easily parallelizable since it samples the posterior for each pixel of interest.}   Suppose that we have multiple vector realizations of $y_j \sim \mathbf{Y} \in \mathbb{R}^n$ to form $Y = [y_1, y_2, y_3, \ldots, y_J]$. In practice, the data vectors $y_j$ are often noisy {indirect} measurements of the same underlying unknown signal represented by a random variable $\mathbf{X}$. It is therefore realistic to suppose that they should share similar support in a chosen sparse domain, for instance in the edge domain.  This is the {\em joint sparsity} assumption.   As noted previously, other sparse domains may also be used for this purpose. 

In order to make use of the joint sparsity concept, we first provide some definitions and notation, which may be found in \cite{VBJSGelb}. Let $p \in \bbR^n$.  As is conventional, we write $\nm{\cdot}_0$ for the $\ell^0$-`norm', i.e.
\bes{
\nm{p}_0 = \left | \mathrm{supp}(p) \right |,
}
where $\supp(p)$ is the support of $p = (p_i)^{n}_{i=1}$ defined by
\bes{
\supp(p) = \{ i : p_i \neq 0 \}.
}
Given a vector $\bm{w} = (w_i)^{n}_{i=1}$ of positive weights, we define the weighted $\ell^1_{\bm{w}}$-norm as
\begin{equation}
    \nm{p}_{1,\bm{w}} = \sum^{n}_{i=1} w_i | p_i |.
\label{eq:weightedl1}
\end{equation}
It is also possible to define weighted $\ell^q_{\bm{w}}$-norms for $q \neq 1$, but this is not needed for our investigation.

If $P = (p_{i,j})^{n,J}_{i,j=1} \in \bbR^{n \times J}$ is a matrix, we define the $\ell^{q_1,q_2}$-norms by
\bes{
\nm{P}_{q_1,q_2} = \left( \sum^{n}_{i=1} \left ( \sum^{J}_{j=1} | p_{i,j} |^{q_1} \right )^{q_2/q_1} \right )^{1/q_2}.
}
Again by convention the $\ell^{2,0}$-`norm' is defined as
\bes{
\nm{P}_{2,0} = \left | \left \{ i : \sum^{J}_{j=1} | p_{i,j} |^2 \neq 0 \right \} \right |.
}
We are now able to provide the following definition: 

\begin{definition}
\label{def:joint_sparsity}
A vector $p \in \bbR^n$ is $s$-sparse for some $1 \leq s \leq N$ if
$$
\nm{p}_{0} = | \supp(p) | \leq s.
$$
A collections of vectors ${p}_1,\ldots,p_J \in \bbR^n$ is $s$-joint sparse if
$$
\nm{P}_{\alpha,0} = \left | \bigcup^{J}_{j=1} \supp \left (p_j \right ) \right | \leq s,\quad\quad{{\alpha > 0}}
$$
where $P = \left [ p_1 | \cdots | p_J \right ]$.  {In typical applications $\alpha$  is chosen to be $2$.}
\end{definition}

Based on our discussion following (\ref{eq:laplace_prior}), the sparsifying operator $\mathcal{L}$ is the PA transform where $p_j = \mathcal{L}x_j$ approximates the solution $x$ in the edge domain for each $j = 1,\cdots,J$. It follows that $P = \left[\mathcal{L}x_1|\cdots| \mathcal{L}x_J\right]$ in Definition \ref{def:joint_sparsity}.   

If the edges were explicitly known, the prior in \eqref{eq:laplace_prior} could be tailored to promote a posterior that not only reflects the existence of the sparse domain but also the regions of its support.  Specifically, a more accurate prior would be dependent on $\supp(\mathcal{L} x)$. Although this information is not directly available in most applications, it is possible to extract it (approximately) from the observable data, as was demonstrated in \cite{VBJSGelb,VBJS, VBJS2} using the VBJS method. In the context of this investigation, the VBJS approach is an example of empirical Bayesian inference, and results in a support informed sparse prior.

\subsubsection{Variance based joint sparsity (VBJS)}\label{subsec:VBJS}
To use the VBJS approach, we first require the calculation of the sample variance $v \in \bbR^n$ across the rows of $P$, where $P  = [\mathcal{L} x_1 | \dots |\mathcal{L}x_J]$ is the sparse matrix defined by Definition \ref{def:joint_sparsity}. The (spatial) entries of $v$ are given by
\begin{equation}
\label{eq:var_comp}
{v}_i = \frac{1}{J}\displaystyle\sum_{j=1}^J {P}_{i,j}^2 - \left(\frac{1}{J}\displaystyle\sum_{j=1}^J {P}_{i,j}\right)^2, \quad i = 1,\dots,n.
\end{equation}
 We note that because we do not have the true solution vector, $x$,  we first must approximate each $\mathcal{L}x_j$ from the observable MMV data $y_j$, $j = 1,\dots,J$.  This may be done in a variety of ways, and in some cases does not require an initial approximation of the solution itself, \cite{VBJS2}.

The VBJS technique constructs a weighting matrix $W = diag(w_1,\dots, w_n)$ that weighs the importance of the prior for different regions of the domain based on the sparse representations of the data. {For the compressive sensing approach used in \cite{VBJSGelb,VBJS}, these spatially varying weights were scaled to essentially be inversely proportional to the reciprocal of \eqref{eq:var_comp}.}  That is,
{
$$
    w_{i} = \frac{1}{v_i + \epsilon},
$$
where $\epsilon > 0$ is some small parameter that offsets zero variance. As discussed in \cite{VBJS}, the weights may be further scaled to account for the different magnitudes of the non-zero values in the sparse domain. While such separation of scales is important for determining the spatially varying weighting parameter in the CS environment, it is not needed in the Bayesian framework for reasons that will become clear below.  What is important is to note that} the VBJS technique's reliance on the variance of the measurements in the sparse domain is motivated by the provable convergence rates of edge detection methods near and away from jumps, \cite{archibald2005polynomial,GT1,GT06,VBJS2}. In particular, if the variance in the sparse domain is high in a region of $\mathbb{R}^n$, then the corresponding entries of $W$ are small. This follows the assumption that high variance values indicate regions where the measurements share support in the sparse domain. Conversely, low variance values suggest truly sparse regions, and in this case the corresponding entries of $W$ should be large to strongly enforce the sparse prior.

For the purpose of empirical Bayesian inference, and specifically to define the support informed sparse prior, we simply threshold the variance so that $W$ becomes a {diagonal} binary mask {$M = diag(m_1,\dots,m_n) \in \mathbb{R}^{n \times n}$} that enforces a prior {in regions that are truly sparse, but assumes no prior in regions of support.}  More specifically, for threshold $\tau = \frac{1}{n}$ we define {the diagonal entries of $M$ as} 
\begin{equation}
    m_{i} = \begin{cases}
    {1}, &\quad w_i \geq \tau \\
    {0}, &\quad w_i < \tau,
    \end{cases}
    \label{eq:mask}
\end{equation}
where $w_i$, $i = 1,\dots, n$, {are provided above.}  
Figure \ref{fig:mask_demo} demonstrates the process for obtaining (\ref{eq:mask}). Observe that because the weights $w_i$ exhibit a large separation of scale, our masking technique is not sensitive to the thresholding parameter $\tau$, {which was chosen simply to be consistent with the resolution of the data}. 

\begin{figure}[htb]
\includegraphics[width = \textwidth]{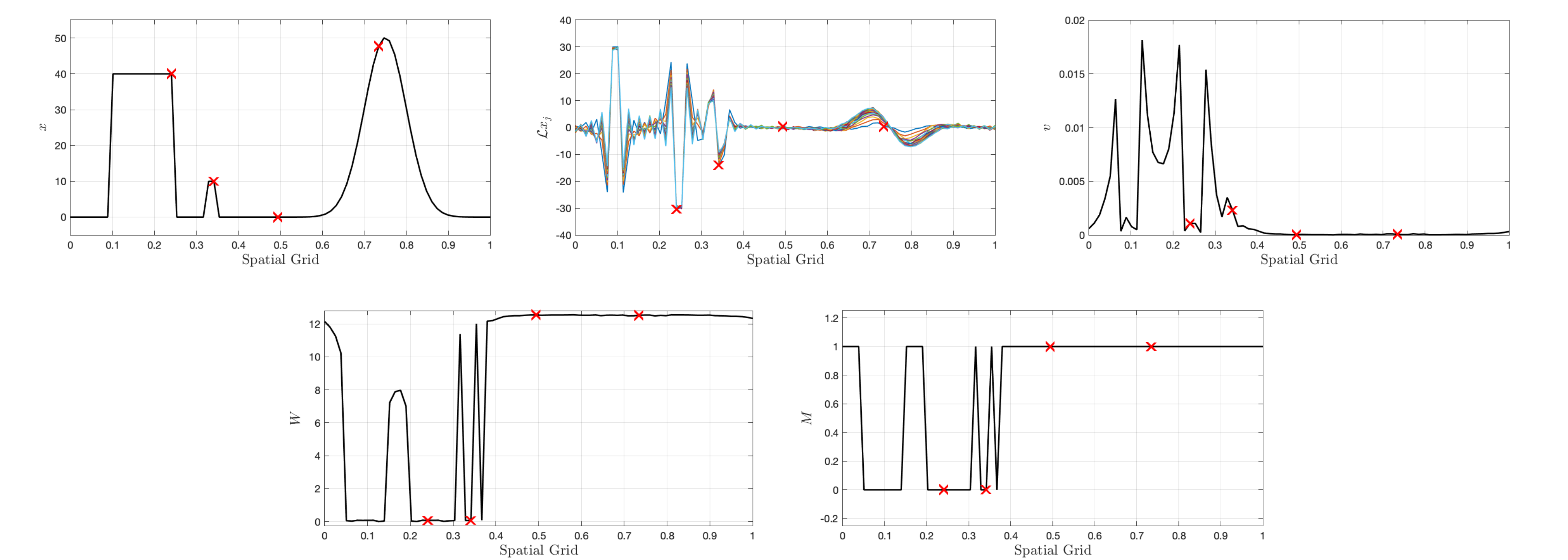}
\caption{Demonstration of the calculation of the mask used in the proposed support informed sparsity prior. (top-left) True function. (top-middle) Jointly sparse vectors generated from indirect, noisy measurements $y_j$, $j = 1,\dots,10$, of $x$. (top-right) Spatial variance $v$ (\ref{eq:var_comp}) across joint sparsity vectors. (bottom-left) Weights calculated using the VBJS technique, \cite{VBJSGelb,VBJS}. (bottom-right) Mask (\ref{eq:mask}) generated for use in the support informed sparsity prior.  The red crosses indicate locations {in the spatial domain $s$} where we calculate evaluation metrics described in Section \ref{sec:numerics} }
\label{fig:mask_demo}
\end{figure}

Given the {binary mask defined} in (\ref{eq:mask}), Algorithm \ref{alg:SISP} describes the VBJS process for constructing a support informed sparsity prior {as it is used for this investigation}. 

\begin{algorithm}[H]
\caption{Support Informed Sparsity Prior}
\label{alg:SISP}
\begin{algorithmic}[1]
\STATE{\textbf{Input:} Set of measurement vectors, ${{ Y}} = [y_1, y_2, \dots, y_J]$.}
\STATE{\textbf{Output:} Support informed sparsity prior.} 
\STATE{Estimate the hyper-parameter $\hat{\lambda}$ using Algorithm \ref{alg:Cross-Validation}.}
\STATE{{Approximate the joint sparsity matrix $P$ defined in Definition \ref{def:joint_sparsity} from the $J$ measurements.}}
\STATE{{Compute the spatial variance $v$ in \eqref{eq:var_comp} from the joint sparsity matrix $P$.}}
\STATE{Use the spatial variance $v$ to determine the weighting vector $w$.}
\STATE{{Use $w$ to determine the binary mask $M$ according to \eqref{eq:mask}.}}
\STATE{Define the support informed sparsity prior density as 
\begin{equation}
    \label{eq:SISP_prior}
    \tilde{f}_{\mathbf{X}|\bm{\hat{\lambda}},M}(x|\hat{\lambda},M) = C \exp\{-{\hat{\lambda} }||{M}\mathcal{L}x||_1\},{\quad C > 0}.%
\end{equation}}
\end{algorithmic}
\end{algorithm}

A couple of remarks are  in order:
\begin{remark}
In the original construction of the weighting matrix $W$ in \cite{VBJS}, the hyper-parameter $\hat{\lambda}$ {was intrinsically fixed}.   Specifically, in the CS solution \eqref{eq:l1reg}, the regularization parameter  $\lambda$ was set to $1$.  {Therefore, to separate the scales of the non-zero components in the sparse domain, it was important that $W$  {\em not} be binary, that is, it had to vary in magnitude to accommodate the various scales of the non-zero components in the sparse domain.}  Furthermore, although not derived this way, one may understand this scaling to be informed by the error estimate of the fidelity term, much like what is obtained through Algorithm \ref{alg:Cross-Validation}.  By contrast, in the current Bayesian framework  we consider the hyper-parameter in the prior to be independent from the support in the sparse domain.  This may be advantageous in circumstances where it is difficult to determine regions of support.  In such cases the resulting mask is close to the identity matrix, and using a hyper-parameter determined through methods such as cross validation may still reduce uncertainty  when compared to standard priors that {have fixed hyper-parameters  (which is analogous to choosing the best regularization parameter in the CS approach).} Hence we set the weighting matrix to be a {binary} mask matrix indicating the support regions and learn the hyper-parameter $\hat{\lambda}$ separately, rather than combining these parameters into one as was done in \cite{VBJS}.  
\end{remark}
{
\begin{remark}
\label{rem:ell2}
The Laplace prior is commonly used to promote sparsity in the sparse domain of the underlying signal, especially when {constructing the MAP estimate of the posterior.} Analogously, $\ell_1$ regularization is chosen in the compressive sensing formulation since using  $\ell_2$ does not effectively promote sparsity, \cite{CompressiveSensingCandes}.   However, when the mask matrix $M$ is applied, the {sparse} regions of $\mathbf{X}$ that are to be regularized are identically zero.  Hence there is no fundamental reason to choose the $\ell_1$ norm.   Indeed, the $\ell_2$ norm was employed in the compressive sensing algorithms designed in \cite{churchill2018edge} and \cite{VBJS}, and in both cases, the resulting algorithms more efficient.  Adopting the same rationale here, for comparative purposes, we also consider a support informed $\ell_2$ prior.
\end{remark}



\subsection{Defining the Posterior} 
\label{sec:post}

Armed with both a likelihood (\ref{eq:LikelihoodDensity}) and the support informed sparsity prior (\ref{eq:SISP_prior}), we are able to define the full posterior (\ref{eq:f_hat}) for the model in (\ref{eq:ForwardModel}).

{Based on our previous discussion, we will consider two posterior estimates which both use \eqref{eq:LikelihoodDensity} and hyper-parameter $\hat{\lambda}$ calculated in Algorithm \ref{alg:Cross-Validation}:} 
\begin{subequations}
\label{eq:final_posterior}
\begin{equation}
\label{eq:PosteriorDensitylaplace}
    \hat{f}_{\mathbf{X}, \bm{\lambda}|\mathbf{Y}}(x, \hat{\lambda}|y)= C_L \exp \{ - {\hat{\lambda}} ||\mathcal{L}x||_1 - \frac{1}{2\sigma^2}||y-Ax||_2^2\},\quad\quad C_L> 0,
\end{equation}
\begin{equation}
\label{eq:PosteriorDensitymask}
    \hat{f}_{\mathbf{X}, \bm{\lambda}|\mathbf{Y}}(x, \hat{\lambda}|y)= C_M \exp \{ - {\hat{\lambda}} ||M\mathcal{L}x||_1 - \frac{1}{2\sigma^2}||y-Ax||_2^2\},\quad\quad C_M > 0.
\end{equation}
\end{subequations}
The first uses the Laplace prior \eqref{eq:laplace_prior}, while the second uses the support informed sparse {$\ell_1$} prior \eqref{eq:SISP_prior} constructed in Algorithm \ref{alg:SISP}.
Although parameters $C_L$ and $C_M$ are not explicitly known, they are not needed in calculation. {Note that within this MMV framework, $\sigma{^2}$ can simply be estimated as the average sample variance across all MMVs in the data space}, {obtained as
\[{\sigma}^2 \approx \frac{1}{n}\left(\frac{1}{J}\displaystyle\sum_{j=1}^J {Y}_{i,j}^2 - \left(\frac{1}{J}\displaystyle\sum_{j=1}^J {Y}_{i,j}\right)^2\right),\]
where ${Y}  = \left [ y_1 | \cdots | y_J \right ]$.
}
For comparative purposes, we will also consider the following posteriors,
\begin{subequations}
\label{eq:final_posteriorl2}
\begin{equation}
\label{eq:PosteriorDensityhierarchical}
    \hat{f}_{\mathbf{X}, \bm{\lambda}|\mathbf{Y}}(x, \hat{\lambda}|y)= C_L \exp \{ - \frac{\hat{\lambda}}{2} ||\mathcal{L}x||_2^2 - \frac{1}{2\sigma^2}||y-Ax||_2^2\},\quad\quad C_L> 0,
\end{equation}
\begin{equation}
\label{eq:PosteriorDensitymaskl2}
    \hat{f}_{\mathbf{X}, \bm{\lambda}|\mathbf{Y}}(x, \hat{\lambda}|y)= C_M \exp \{ - \frac{\hat{\lambda}}{2} ||M\mathcal{L}x||_2^2 - \frac{1}{2\sigma^2}||y-Ax||_2^2\},\quad\quad C_M > 0,
\end{equation}
\end{subequations}
{for which we use the same approximation of $\sigma^2$ as above.} As discussed in the introduction, the method in \cite{BardsleyUQ} uses a hierarchical Gaussian prior (HGP) with the scalar hyper-priors chosen to be Gamma distributed. Furthermore, the sparse transform operator in \cite{BardsleyUQ}  is symmetric positive definite.  Because of the resulting conjugacy relationships, the posteriors are known in closed form and the Gibbs sampling method can be efficiently applied, which is a primary objective in \cite{BardsleyUQ}.  {It is important to note that employing the mask matrix $M$ to ${\mathcal L}$ (even if ${\mathcal L}$ is constructed to be SPD),  would not yield an SPD transform.}  Without considering {the convergence properties of the corresponding numerical implementations}, the posterior provided in \eqref{eq:PosteriorDensityhierarchical} is comparable to the one in \cite{BardsleyUQ}.  Indeed it is suggested there that additional information {may be used to determine the hyper-prior.  The $K$-fold  cross validation algorithm in  Algorithm \ref{alg:Cross-Validation} provides such information.} On the other hand}, {when the underlying signal is itself sparse so that ${\mathcal L} = \mathcal{I}$ the identity matrix,} \eqref{eq:PosteriorDensitymaskl2}, which uses the support informed  $\ell_2$ prior, can be viewed as  alternative approach to SBL (using EM), \cite{tipping2001sparse, wipf2004sparse}, since both our new algorithm and SBL {incorporate the support locations} in the sparse domain to inform the prior.  Additional contrasts and comparisons are provided in the introduction.
Finally we note that while using Gibbs sampling is generally more efficient, it is not consistent with employing the support informed prior which reduces the uncertainty of the solution.  Moreover, our method will allow downstream processing to be more efficient, as we can more efficiently probe regions of interest.

{Finally, we point out that we can incorporate the MMV framework directly into the likelihood estimate that appears in each of the posterior estimates, \eqref{eq:final_posterior} and \eqref{eq:final_posteriorl2}.  However, as we are only considering one data source given by (\ref{eq:ForwardModel}), this information is somewhat redundant.  Furthermore, the added MMV information is already accounted for in determining the support of the underlying signal.  In our numerical experiments we use $y = \frac{1}{J}\sum_1^J y_j$, i.e. the mean of the observations,  in each posterior estimate.}

%% file: posterior.tex
{Obtaining a  maximum a posteriori (MAP) estimate of (\ref{eq:final_posterior}) or \eqref{eq:final_posteriorl2} is relatively straightforward. However, because it provides only a mode of the desired distribution, the MAP estimate is unable to quantify the uncertainty of the results. 
Moreover, a mode of the posterior distribution may be an atypical solution in the sense that {even if it is unique,} it may differ substantially from its mean or median. 
By contrast, the sampling techniques described in Section \ref{sec:IntervalEstimation} allow for both uncertainty quantification in the signal recovery as well as the ability to discern modes of the posterior distribution.}

\subsection{Point Estimation}
\label{sec:PointEstimation}
\label{sec:map}

{While the main goal of this investigation is to draw samples of the full posterior, it is instructive to first demonstrate how the MAP point estimate is used to recover a mode of the unknown distribution,  \cite{murphy2012machine}.}  

For any realization of \eqref{eq:ForwardModel},  $y_j \sim \mathbf{Y}$, $j  = 1,\dots,J$, the MAP estimate for either posterior in  \eqref{eq:final_posterior} or \eqref{eq:final_posteriorl2} is computed as
\begin{equation}
\label{eq:MAP1}
    x_{MAP}^j \equiv \argmax{x \in \mathbb{R}^n} \left\{\log \left[ \hat{f}_{\mathbf{X}, \bm{\lambda}|\mathbf{Y}}(x, \hat{\lambda}|y_j)\right]\right\}, \quad j = 1,\dots,J.
    \end{equation}
A corresponding convex optimization problem is then easily derived for the posteriors in \eqref{eq:final_posterior} or \eqref{eq:final_posteriorl2} as 
\begin{align}
x_{MAP}^j &=
   \argmax{x \in \mathbb{R}^n} \Big\{ \log(C) - \Big( \frac{\hat{\lambda}}{\nu} ||\mathcal{T}x||_\nu^\nu + \frac{1}{2 \sigma^2}||y_j-Ax||_2^2 \Big) \Big\} \nonumber \\
    &= \argmin{x \in \mathbb{R}^n} \Big\{ \frac{\hat{\lambda}}{\nu} ||\mathcal{T}x||_\nu^\nu + \frac{1}{2 \sigma^2} ||y_j-Ax||_2^2 \Big\}, \quad j = 1,\dots,J, \label{eq:MAPEstimate}
\end{align}
{where $\nu = 1$ or $2$, $\mathcal{T} = \mathcal{L}$ or $\mathcal{ML}$, and $C = C_L$ or $C_M$.}
 Since the MAP estimation can be reduced to an optimization problem, \eqref{eq:MAPEstimate} is often understood from a non-Bayesian perspective.  {In particular, note that the MAP estimate used in Step 7 of Algorithm \ref{alg:Cross-Validation} is similarly constructed, with $\hat{\lambda}$ given by $\tilde{\lambda}$ in Step 6, $\mathcal{T} = \mathcal{L}$, and $\nu = 1$.   We will also use the average MAP estimate given by  (\ref{eq:MAPEstimate}) over the $J$ measurement vectors to set the initial state $x^0$ for the Metropolis-Hastings Algorithm given in Algorithm \ref{alg:MHAlgorithm}.  Finally, although not discussed in this current investigation, we point out that there may be situations for which there are regions in a signal domain where efficiency is more important than quantifying uncertainty.  In such cases using \eqref{eq:MAPEstimate} to construct an approximation within these regions may be more appropriate.}


\subsection{{Drawing} samples of the posterior distribution}
\label{sec:IntervalEstimation}
We seek to sample the posterior probability distribution of $\mathbf{X}$ described by (\ref{eq:final_posterior}) or \eqref{eq:final_posteriorl2}  for the model given by \eqref{eq:ForwardModel}. This allows us both to quantify the uncertainty in our signal recovery as well as to discern modes of the posterior distribution that are unavailable in a MAP estimate.  By Bayes' rule, \eqref{eq:full_posterior}, it is clear that {sampling} $\mathbf{X}$ in \eqref{eq:ForwardModel} from observations in \eqref{eq:PointEstimate} requires an approximation to the probability density $f_{\mathbf{Y}}(y)$,  which typically does not have a closed form solution. Moreover, numerical quadrature rules are impractical whenever the dimension of the sample space $\mathbb{R}^n$ is high \cite{kaipio2006statistical}. We therefore choose to use the Markov Chain Monte Carlo (MCMC) method which is designed to construct an ergodic Markov chain that converges to the stationary probability density ${f}_{\mathbf{X},\bm{\lambda}|\mathbf{Y}}(x,\lambda|y)$ in \eqref{eq:final_posterior}, which is proportional to the posterior we seek. Overviews of MCMC may be found in \cite{MCMC,MCMCR, kaipio2006statistical}. Because our {approach}  does not assume conjugate priors,  the posterior density (\ref{eq:full_posterior}) we derive does not represent a known posterior distribution. Hence we employ the particular MCMC method known as the Metropolis-Hastings (MH) algorithm \cite{hastings, metropolis}, {which we briefly describe below}.

\subsubsection{The Metropolis Hastings (MH) Algorithm}
\label{sec:MH}

In the MH algorithm, a Markov chain in ${\Omega}$ is constructed so that it converges to the conditional posterior density ${f}_{\mathbf{X},\bm{\lambda}|\mathbf{Y}}(x,\lambda|y)$  in \eqref{eq:full_posterior}. 
The MH algorithm first defines a \emph{proposal distribution} $q: \Omega \times \Omega \rightarrow [0, 1]$ such that when the Markov chain is at state $x^{k-1}$, it proposes a \emph{candidate} state $x^{cand} \sim q$. Its acceptance of the next state, $x^k$, is determined by an \emph{acceptance probability} $\alpha: \Omega \times \Omega \rightarrow [0,1]$ given by

\begin{align}
    \alpha(x^{cand}|x^{k-1}) 
    &= \min \Bigg\{ 1, \frac{q(x^{k-1}|x^{cand}){f}_{\mathbf{X},\bm{\lambda}|\mathbf{Y}}(x^{cand},\lambda|y)}{q(x^{cand}|x^{k-1}){f}_{\mathbf{X},\bm{\lambda}|\mathbf{Y}}(x^{k-1},\lambda|y)} \Bigg\}\nonumber\\
    &{\approx} \min \Bigg\{ 1, \frac{q(x^{k-1}|x^{cand})\hat{f}_{\mathbf{X},\bm{{{\lambda}}}|\mathbf{Y}}(x^{cand},{\hat{\lambda}}|y)}{q(x^{cand}|x^{k-1})\hat{f}_{\mathbf{X},\bm{{{\lambda}}}|\mathbf{Y}}(x^{k-1},{\hat{\lambda}}|y)} \Bigg\},
\label{eq:acceptance}
\end{align}
where  {$\hat{f}_{\mathbf{X},{\bm{\lambda}}|\mathbf{Y}}$ is defined in \eqref{eq:final_posterior}}.  In this way the MH algorithm is able to evaluate the full posterior (\ref{eq:full_posterior}) without needing to calculate the evidence $f_{\mathbf{Y}}(y)$ in (\ref{eq:MarginalDensity}). The general MH algorithm is summarized in  Algorithm \ref{alg:MHAlgorithm}.

\begin{algorithm}[H]
\caption{The Metropolis-Hastings Algorithm}
\label{alg:MHAlgorithm}
\begin{algorithmic}[1]
\STATE{\textbf{Input:} The initial state $x^0$, {the choice of proposal distribution $q$, and its corresponding variance $\sigma_q$.}}
\STATE{\textbf{Output:} Markov chain $\{x^1, x^2, x^3, \ldots, x^M\}$.}
\STATE{Initialize $x^0$ and $k = 1$.}
\WHILE{$k\leq M$}
	\STATE{Propose $x^{cand} \sim q(x^k|x^{k-1})$.}
 \STATE{Compute $\alpha(x^{cand}|x^{k-1})$ using (\ref{eq:acceptance}) or (\ref{symmetricproposal})}.
    \STATE{Sample $u \sim \mathcal{U}[0, 1]$.}
    \IF{$u < \alpha$}
        \STATE{Accept the candidate: set $x^k = x^{cand}$.}
	\ELSE
        \STATE{Reject the candidate:  set $x^k = x^{k-1}$.}
   \ENDIF
\ENDWHILE
\end{algorithmic}
\end{algorithm}

Algorithm \ref{alg:MHAlgorithm} requires user inputs of the initial state $x^0$  and the proposal distribution $q$. For our experiments we set the initial state $x^0$ as the average MAP estimate in (\ref{eq:MAPEstimate}) over all measurement vectors $y_j$, $j = 1,\dots, J$, and $q$ to be Gaussian, such that $x^{cand} \sim \mathcal{N}(x^{k-1}, \sigma_q)$. {As noted previously, we also use $y = \frac{1}{J} \sum_{j = 1}^J y_j$, the mean of the observable data.}
Since the Gaussian distribution is symmetric, i.e.~$q(x^{k-1}|x^{cand}) = q(x^{cand}|x^{k-1})$,  the acceptance probability (\ref{eq:acceptance}) is simplified to 
\begin{equation}\label{symmetricproposal}
    \alpha(x^{cand}|x^{k-1}) = \min \Bigg\{ 1, \frac{\hat{f}_{\mathbf{X},\bm{\lambda}|\mathbf{Y}}(x^{cand},\hat{\lambda}|y)}{\hat{f}_{\mathbf{X},\bm{\lambda}|\mathbf{Y}}(x^{k-1},\hat{\lambda}|y)} \Bigg\}.
\end{equation}
While the ideal acceptance ratio is an open question for the Metropolis-Hastings algorithm, a typical range for an acceptance ratio is $[.2,.8]$. This enables the method to ``reasonably'' explore  the space, preventing the chain from getting stuck in one place or moving too quickly from one component to another, \cite{MCMCR}.  Because the acceptance ratio explicitly depends on $\sigma_q$, for the purposes of this investigation, we simply chose $\sigma_q$ so that it fell within this interval.  Specifically, for all posteriors given in (\ref{eq:final_posterior}) and \eqref{eq:final_posteriorl2}, we ran a few iterations of Algorithm \ref{alg:MHAlgorithm} and observed the behavior of  the acceptance ratio.  If it fell outside the interval, we adjusted the input parameter $\sigma_q$ accordingly.  

It is important to note that $x^0$ may not provide a good initial guess.  This is a standard problem in all MCMC algorithms and is addressed by implementing what is known a {\em burn-in} period.  That is, if $M$ is the total number of iterations computed by Algorithm \ref{alg:MHAlgorithm}, the first $B< M$ (called the burn-in rate) iterations  are discarded, so that the total length of the chain is $M-B$.  {We chose $M = 50,000$ and $B = 25,000$ in all of our experiments, the same values used in studies in which convergence is consistently  observed, \cite{ConvergenceMCMC}.}  
In this regard we point out that our strategy for choosing $x^0$, and more importantly our use of the support informed sparse prior, allowed us to pick smaller burn-in rate and chain length while still reducing the uncertainty in all of our experiments.  More studies are needed to confirm this is the case in general.  
%


The iterations of the Markov chain obtained using Algorithm \ref{alg:MHAlgorithm} may be represented by a matrix $\mathcal{M} \in \mathbb{R}^{n \times M}$ {given by}
\begin{equation}
\label{eq:MarkovMatrix}
\mathcal{M}=
\begin{bmatrix}
    x_{1}^0 & x_{1}^1 & x_{1}^2 & \dots  & x_{1}^M \\
    x_{2}^0 & x_{2}^1 & \dots  & x_{2}^3 \\
    \vdots & \vdots & \vdots & \ddots & \vdots \\
    x_{n}^0 & x_{n}^1 & x_{n}^2 & \dots  & x_{n}^M
\end{bmatrix},
\end{equation}
where each column of the matrix $\mathcal{M}$ represents one state in the Markov chain. In particular, an entry $x^t_i$ is the $i^{th}$ {component} of the $t^{th}$ state in the Markov chain. As in \cite{BardsleyUQ},  we then calculate a final estimate for the unknown $x$ {at each component} as the mean of the samples after discarding the first $B$ states in the Markov chain:
\begin{equation}\label{eq:final_approx}
    x_i= \frac{1}{M -B}\sum_{t = B+1}^M x_i^t.
\end{equation}

A summary of our proposed empirical Bayesian inference approach to solving the statistical inverse problem defined in (\ref{eq:ForwardModel}) is provided in Algorithm \ref{alg:final_alg}.  {The mean in \eqref{eq:final_approx} is obtained after discarding the first $B < M$ burn-in iterations.} {The Algorithm is written for the support informed sparse $\ell_1$ prior, but the $\ell_2$ prior case given in \eqref{eq:PosteriorDensitymaskl2} follows analogously.}

\begin{algorithm}[H]
\caption{{Support Informed Sparse Prior} Empirical Bayesian Inference Algorithm}
\label{alg:final_alg}
\begin{algorithmic}[1]
\STATE{\textbf{Input:} Multiple measurement vectors, $y_j$, for $j = 1,...,J$.}
\STATE{\textbf{Output:} Markov chain $\mathcal{M}$ in \eqref{eq:MarkovMatrix} describing samples from $f_{\textbf{X},\bm{\lambda}|\mathbf{Y}}(x,\lambda|y)$.}
\STATE{Estimate the hyper-prior $\hat{\lambda}$ according to Algorithm \ref{alg:Cross-Validation}.}
\STATE{Calculate the support informed sparsity prior using Algorithm \ref{alg:SISP}.}
\STATE{Sample the posterior (\ref{eq:PosteriorDensitymask}) using Algorithm  \ref{alg:MHAlgorithm}.}
\end{algorithmic}
\end{algorithm}

%

%% file: numerics.tex
We are now ready to demonstrate the utility of our support informed sparsity prior (\ref{eq:SISP_prior}) in sampling the full posterior (\ref{eq:PosteriorDensitymask}) {and the posterior derived using the support informed $\ell_2$ prior used in \eqref{eq:PosteriorDensitymaskl2}.} 
{In analyzing our results, we will calculate the relative error
\begin{equation}
\label{eq:rel_error}
\frac{||x-x^{true}||}{||x^{true}||}.
\end{equation}
Here $x$ is determined by \eqref{eq:final_approx} and $x^{true}$ is the underlying signal.  We will also consider measurements that are acquired with increasing levels of noise.  The  signal to noise ratio (SNR) computed in decibals (dB) according to
\begin{equation}
\label{eq:SNR}
\text{SNR} = 10\log_{10}\left(\frac{\mathbb{E}[\mathbf{Y}]}{\sigma^2}\right),
\end{equation}
where $\mathbb{E[\mathbf{Y}]}$ is the expected value of the collected data and $\sigma^2$ is the variance of the complex additive Gaussian white noise (AGWN) in (\ref{eq:ForwardModel}).}

To quantify the uncertainty of our recovery,  we  approximate the credibility interval (also called the Bayesian interval) {of the resulting distribution}, \cite{murphy2012machine, kaipio2006statistical, casella2002statistical}, which is defined for {$x \in \mathbb{R}$} as:
\begin{definition}\label{def:credint}
Given $\eta \in (0,1)$, the credibility set $C_{\eta} = (\ell,u)$ is an equiprobability hypersurface that contains $100(1-\eta)$\% of the mass of the posterior distribution and is defined according to 
$$\int_{C_\eta} \hat{f}_{\mathbf{X},\bm{\lambda}} (x,\hat{\lambda}|y) dx = 1-\eta, \quad \hat{f}_{\mathbf{X},\bm{\lambda}} (x,\hat{\lambda}|y)|_{x\in \partial C_\eta} = \text{constant}.$$ 
\end{definition}
{Following Definition \ref{def:credint}, we calculate $C_\eta$ as}
\begin{equation}\label{eq:CredibilityInterval}
    C_\eta = \{(\ell ,u) \hspace{1mm} | \hspace{1mm} \hat{f}_{\mathbf{X},\bm{\lambda} | \mathbf{Y}}(\ell \le x \le u,\hat{\lambda} | y) = 1-\eta\}.
\end{equation}
For example, {$C_{0.05}$} is a $95\%$  credibility interval of a probability density and is interpreted as a sub-interval of the domain that contains $95\%$  of the probability density. 
We note that $C_\eta$ is generally not unique because there is typically more than one interval that satisfies \eqref{eq:CredibilityInterval}. This fact emphasizes the importance of initiating the MCMC chain with an appropriate point estimate, so that we may constrain our recovered credibility interval to contain the mode, \cite{murphy2012machine}.
{Hence, as noted previously, we choose} the beginning state of the chain $x^0$ in Algorithm \ref{alg:MHAlgorithm} to be the average MAP estimate \eqref{eq:MAPEstimate} over all measurements $y_j$ for $j = 1,\dots,J$, which is approximately the peak of the mode of the posterior probability density. In this way we ensure that the Markov chain explores a meaningful region of the posterior density, \cite{MCMCR}. 




We also analyze the convergence of our MCMC chains to steady state through the use of autocorrelations, trace plots and acceptance ratios. 
Trace plots show the state of the chain over each iteration. 
Theoretically, if the Markov chain has converged to the posterior density, the states of the chain should not be correlated. {Therefore, the chain is said to be {\em well-mixed} if the  trace plot exhibits behavior indicating that the states are independently and identically distributed (iid).   On the other hand,} any portion of the Markov chain having high auto-correlation between states would imply that the Markov chain has not {(yet)} converged to the stationary distribution.  {As another indicator of mixing}, we use the  Autocorrelation Function (ACF) package in Matlab, \cite{acf_matlab}, to compute the correlogram of each Markov Chain recovered through the proposed MCMC sampling techniques.

The acceptance ratio {determined from line 9} in Algorithm \ref{alg:MHAlgorithm} is the number of times $x^{cand}$ is accepted as the next sample relative to the total number of iterations.    Therefore, a large acceptance ratio {(typically $\alpha > .8$)} may indicate the Markov chain is stuck in a certain region of the posterior probability density and has not explored the entire support of the density. In contrast, a small acceptance ratio {(e.g.~$\alpha < .2$)} may indicate that the chain may be moving too quickly within the domain of the {proposal} density.   Further, with a low acceptance ratio, it remains possible that the chain does not explore certain isolated modes of the posterior probability density \cite{MCMCR}.   While the ideal acceptance ratio remains an open question for the Metropolis-Hastings algorithm, as a general rule of thumb, it should be close to $.5$ for a one-dimensional problem \cite{MCMCR}. This metric was primarily derived for the Gaussian proposal distribution, which is also used in this paper.

\subsection{Sparse Signal Recovery in One Dimension}
\input{sparsesignal}

\subsection{Signal Recovery in One Dimension}
We now consider reconstructing the signal defined as 
\begin{example}\label{ex:true_1D}
\begin{equation*}
    x(s) = \begin{cases}
    40, &\quad 0.1 \leq s \leq 0.25 \\
    10, &\quad 0.35 \leq s \leq 0.325 \\
    20\sqrt{2 \pi}e^{-\left(\frac{s-0.75}{0.05}\right)^2}, &\quad s > 0.5 \\
    0,  & \quad \mathrm{otherwise,}
    \end{cases}
\end{equation*}
\end{example}
where $s \in [0,1]$ is the domain of interest. Similar to {the example used in} \cite{BardsleyUQ}, the domain is uniformly discretized so that $s_i = \frac{i}{n}$, $i = 1,\dots,n$ with $n = 80$. For all of our experiments in one dimension, we simulate $J=20$ measurements. {We choose $20$ since it was shown to be a sufficient number of measurements for reconstructing a signal using the VBJS method in \cite{VBJS}.}  We also choose  $K=20$ and $M=10$  for Algorithm \ref{alg:Cross-Validation}.

Because we are {ultimately} interested in imaging applications such as MRI, SAR and ultrasound, in our experiments we assume that the observable data {$y$} consists of multiple, noisy Fourier measurements of $x$. That is, each measurement is collected according to (\ref{eq:ForwardModel}), with $A$ defined as the discrete Fourier transform matrix and $\sigma >0$. {We point out that}  the discrete Fourier transform matrix is better conditioned than the blurring model used in \cite{BardsleyUQ}.  Here we allow for the measurements to be collected with significantly higher noise levels, i.e., $\sigma >>0$ {or SNR $<<1$}, however.  {In what follows we show how results for $SNR = 5$ dB and $SNR \approx -.2$ dB.}

\begin{figure}[htb]
\centering
\includegraphics[width=\textwidth]{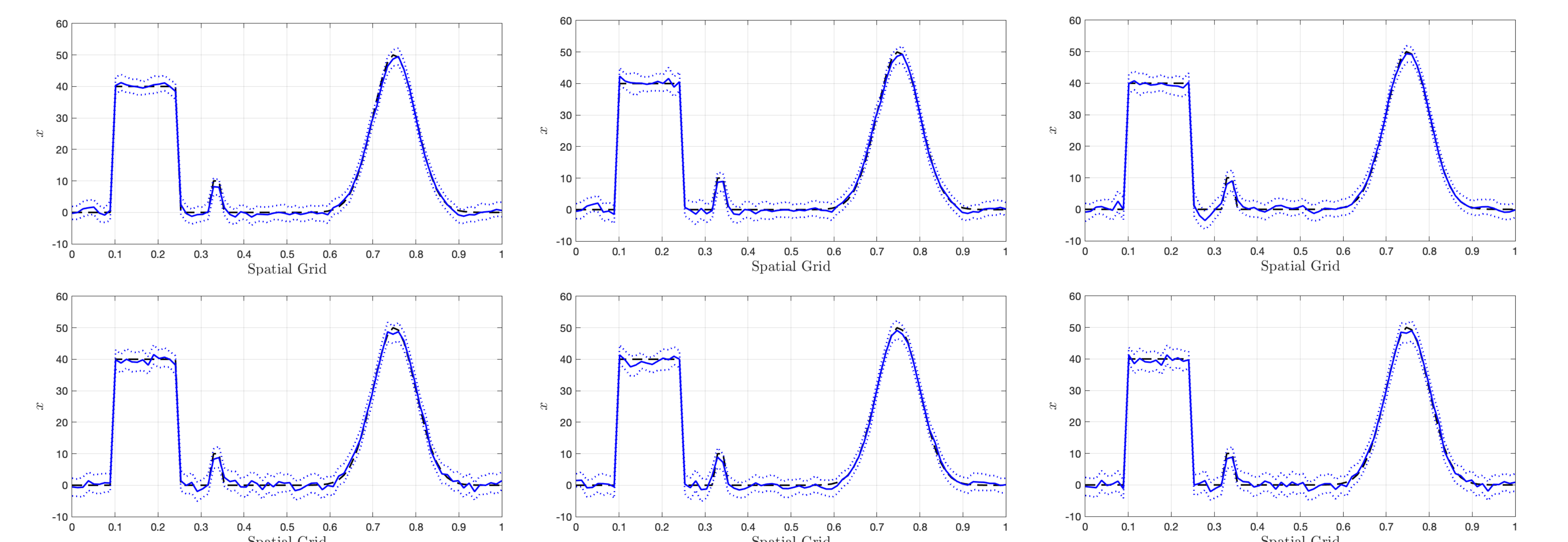}
\caption{Approximate reconstructions (\ref{eq:final_approx}) of Example \ref{ex:true_1D} displayed along with credibility intervals $C_{0.05}$ when using (left) the Laplace prior (\ref{eq:PosteriorDensitylaplace})  (middle) the proposed support informed sparsity prior (\ref{eq:PosteriorDensitymask}) and (right) the proposed support informed prior described in (\ref{eq:PosteriorDensitymaskl2}).  {(top) $\sigma \approx 5.5$ ($SNR = 5$ dB), (bottom) $\sigma = 10$ ($SNR \approx -.2$ dB).}}
\label{fig:recons}
\end{figure}

Figure \ref{fig:recons} displays the mean of the Markov chain (\ref{eq:final_approx}) recovered on each grid point in the domain when sampling the posteriors described in (\ref{eq:final_posterior}) with (left) the Laplace prior (\ref{eq:laplace_prior}) and (right) the proposed support informed sparsity prior (\ref{eq:SISP_prior}) for {(top) $\sigma \approx 5.5$ ($SNR = 5$ dB) and (bottom) $\sigma = 10$ ($SNR \approx -.2$ dB)}. Along with the mean, in Figure \ref{fig:recons} we also display $C_{0.05}$ for each reconstruction.  Notice that the proposed support informed sparsity prior results in tighter credibility intervals and more accurate reconstructions, {especially in the $\sigma = 10$ case}.

\begin{figure}[htb]
\centering
\includegraphics[width=\textwidth]{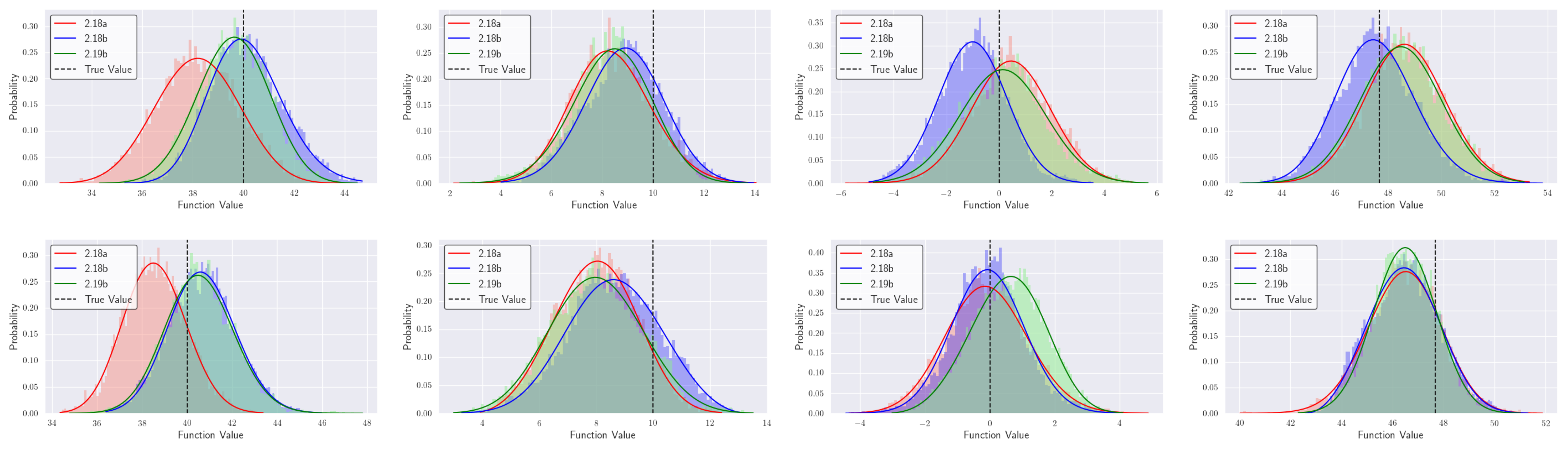}
\caption{{Estimated posterior distributions at the points in the domain indicated by the red crosses in Figure \ref{fig:mask_demo}. Table \ref{tab:distros} reports each distribution and its corresponding parameters. (top) $\sigma \approx 5.5$ ($SNR = 5$ dB); (bottom) $\sigma = 10$ ($SNR \approx -.2$ dB).}}
\label{fig:distros}
\end{figure}

Figure \ref{fig:distros} shows the resulting distributions at the points in the domain indicated by the red crosses in Figure \ref{fig:mask_demo}. Each distribution is estimated given all samples of the Markov chain at the desired point in the domain after discarding the first $B$ states. Distributions are optimized based on the Pearson Chi-Squared goodness of fit metric and the mean-squared error when compared to the histogram. The resulting distributions and their corresponding parameters are displayed in Table \ref{tab:distros}. Consistent with the remainder of our results, we see that the distributions of the support informed $\ell_1$ prior at edge locations is more closely centered around the true function values. In smooth regions, we see compelling evidence to use the support informed $\ell_2$ prior. This motivates future work into exploring the use of different priors in different regions of the domain.

\begin{table}[hbt]
\begin{center}
\resizebox{\textwidth}{!}{
\begin{tabular}{|c|c|c|c|c|c|}
\hline
Posterior & Noise & Location 1 & Location 2 & Location 3 & Location 4 \\ \hline 
\multirow{2}{*}{(\ref{eq:PosteriorDensitylaplace})} &  \multirow{2}{*}{ 5 dB} & Beta & $\chi^2$ & Normal & Log-Normal \\ 
& & (12.45,11.94) & (247.60) & (0.46, 1.50) & (0.01) \\ \hline
\multirow{2}{*}{(\ref{eq:PosteriorDensitymask})} & \multirow{2}{*}{ 5 dB} & $\chi^2$ & Pearson type III & Normal & Log-Normal \\ 
& & (95.34) & (-0.07) & (-1,1.29) & (0.06) \\ \hline
\multirow{2}{*}{(\ref{eq:PosteriorDensitymaskl2})} & \multirow{2}{*}{ 5 dB} & Beta & Beta & Normal & Normal \\ 
& & (22.15,18.08) & ($5.13\times 10^6$, 196.99) & (0.15, 1.61) & (48.48, 1.53) \\ \hline 
\multirow{2}{*}{(\ref{eq:PosteriorDensitylaplace})} &  \multirow{2}{*}{ -.2 dB} & Beta &Weibull (min) & Beta & Pearson type III \\ 
& & (23.53,34.94) & (3.61) & (399.66,$1.15\times 10^7$) & (-0.11) \\ \hline 
\multirow{2}{*}{(\ref{eq:PosteriorDensitymask})} & \multirow{2}{*}{ -.2 dB} & Log-Normal & Weibull (max) & Log-Normal & Weibull (min) \\ 
& & (0.03) & (3.32) & (0.01) & (3.2) \\ \hline
\multirow{2}{*}{(\ref{eq:PosteriorDensitymaskl2})} & \multirow{2}{*}{ -.2 dB} & Beta & Beta & Weibull (max) & Beta \\ 
& & (325.45,$3.32\times 10^3$) & (16.56,19.12) & (3.9) & (48.37,44.73) \\ \hline 
\end{tabular}}
\end{center}
\caption{{The distributions and corresponding parameters resulting from a best-fit analysis of the last $M-B$ Markov chain samples (\ref{eq:MarkovMatrix}) at the locations in the domain indicated by the red crosses seen in Figure \ref{fig:mask_demo}. The distributions along with the histogram of samples are visualized in Figure \ref{fig:distros}.}  }
\label{tab:distros}
\end{table} 

\begin{remark}
We do not report on results using the posterior given by \eqref{eq:PosteriorDensityhierarchical} for Example \ref{ex:true_1D} because the MH algorithm {as implemented in Algorithm \ref{alg:MHAlgorithm}}  did not converge in this case.  We speculate that the (unmasked) $\ell_2$ norm in the prior is prohibitively large for numerical computation. {It may be possible to calculate \eqref{eq:acceptance} so as to avoid this issue, for example by using logarithms.}  Alternatively, one could use a Gibbs sampler with an appropriate SPD operator to replace ${\mathcal L}$ in \eqref{eq:PosteriorDensityhierarchical}, see for example the discussion in \cite{BardsleyUQ}.  {Since our previous observations in the sparse signal case suggest that the resulting confidence interval would still be larger than those displayed in Figure \ref{fig:recons}, and as comparing different sampling algorithms is not the focus of this investigation, we give no further consideration to \eqref{eq:PosteriorDensityhierarchical}.}
\end{remark}

\begin{figure}[htb]
\centering
\includegraphics[width=.9\textwidth]{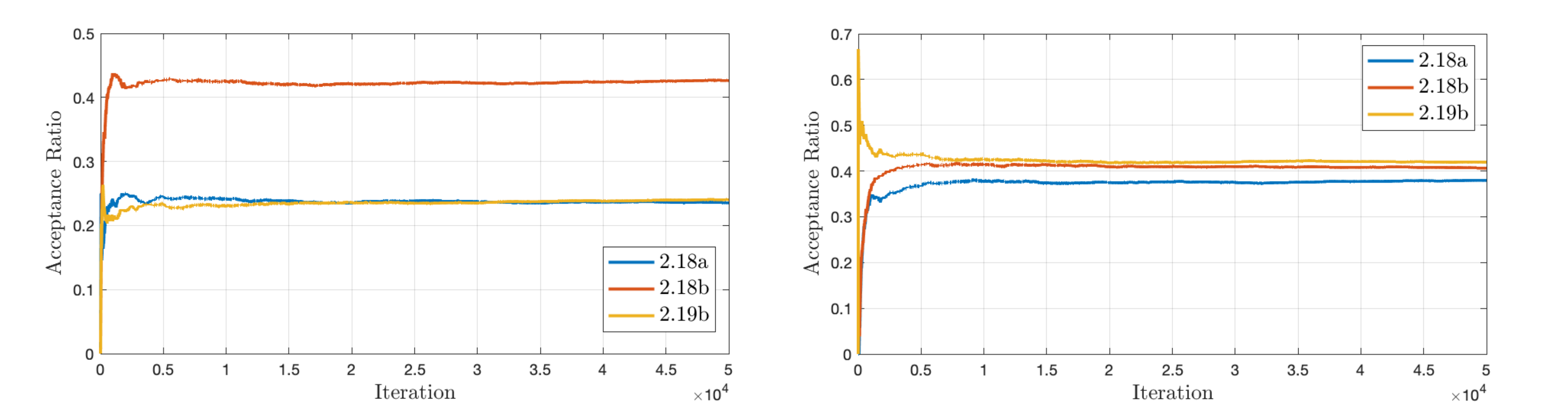}
\caption{Acceptance ratios resulting from sampling the posterior defined using the Laplace prior (\ref{eq:PosteriorDensitylaplace}), support informed sparsity prior (\ref{eq:PosteriorDensitymask}) and support informed $\ell_2$ prior (\ref{eq:PosteriorDensitymaskl2}).  (left) $\sigma \approx 5.5$ ($SNR = 5$ dB); (right) $\sigma = 10$ ($SNR \approx -.2$ dB).}
\label{fig:ar}
\end{figure}

\begin{figure}[htb]
\centering
\includegraphics[width=\textwidth]{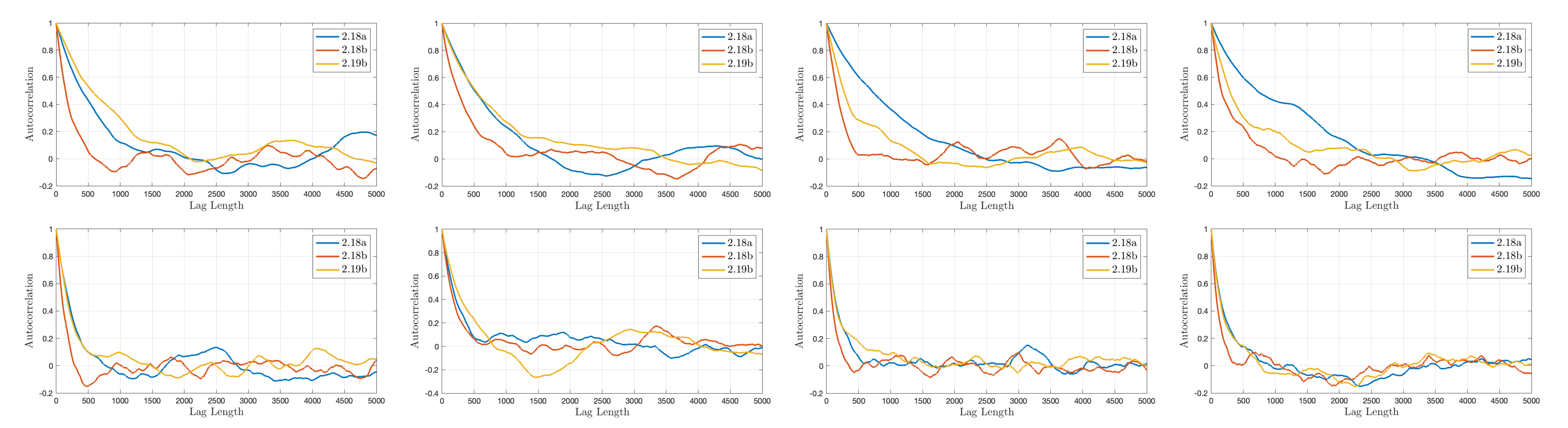}
\caption{The autocorrelation function (ACF) calculated for various lag lengths at the points in the domain indicated by the red crosses in Figure \ref{fig:mask_demo}. Samples are taken from the posterior distributions defined using the  Laplace prior (\ref{eq:PosteriorDensitylaplace}), support informed sparsity prior (\ref{eq:PosteriorDensitymask}) and support informed $\ell_2$ prior (\ref{eq:PosteriorDensitymaskl2}). {(top) $\sigma \approx 5.5$ ($SNR = 5$ dB), (bottom) $\sigma = 10$ ($SNR \approx -.2$ dB).}}
\label{fig:acf_lines}
\end{figure}

\begin{figure}[htb]
\centering
\includegraphics[width=\textwidth]{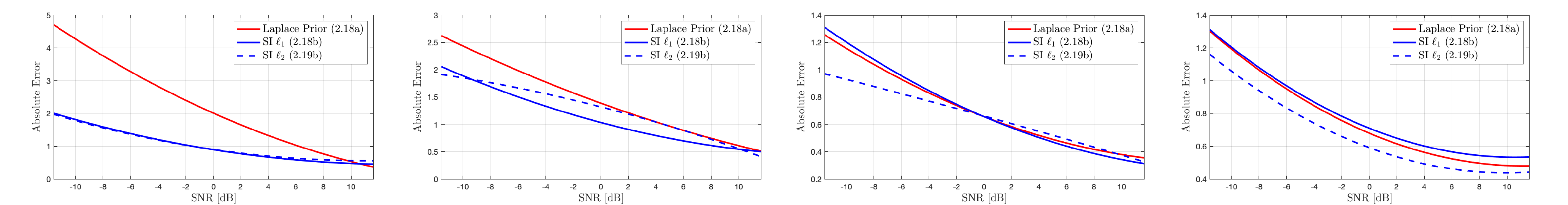}
\caption{Absolute error comparison for decreasing SNR (\ref{eq:SNR}) at each of the red crosses in Figure \ref{fig:mask_demo}}
\label{fig:abs_error}
\end{figure}

We now evaluate the mixing and convergence of the MCMC chains. Figure \ref{fig:ar} displays the acceptance ratio across all MH iterations.  With our strategy to adjust the variance of the proposal distribution, we see that each chain exhibits an acceptable acceptance ratio, with the acceptance ratio for the support informed sparsity prior tending closer to the desirable $.5$ value.   We note that we did not try to optimize $\sigma_q$ in Algorithm \ref{alg:MHAlgorithm}, and that its initial tuning as described in the paragraph below Algorithm \ref{alg:MHAlgorithm} was done in the same way for each choice of posterior.}  Figure \ref{fig:acf_lines} displays the correlograms associated with sampling the posterior distributions at each grid point marked with a red cross in Figure \ref{fig:mask_demo}.  {The results affirm that prior information about the support locations in the sparse domain, which is achievable $SNR = 5$ dB, does indeed reduce the uncertainty.}  When comparing the correlograms resulting from {the posteriors in \eqref{eq:PosteriorDensitylaplace}, \eqref{eq:PosteriorDensitymask}, and \eqref{eq:PosteriorDensitymaskl2} for lower SNR, it is evident that the proposed support informed priors yield less or similarly correlated samples near discontinuities than the more standard Laplace prior does.} Consistent with the results in the sparse signal case, the Laplace prior demonstrates earlier mixing in smooth regions than the support informed sparse prior, but {the mixing is comparable to that achieved using} the support informed $\ell_2$ prior.  {These results, confirmed in Figure \ref{fig:abs_error},}  suggest that as the SNR decreases a better result overall may be obtained when the choice of posterior is itself domain dependent; that is, when the support informed $\ell_2$ prior and the support informed sparse prior are adopted in different regions once the support in the sparse domain has been approximated.  This will be discussed in future investigations, as such a ``splitting'' may also yield more efficient computation.

\begin{figure}[htb]
\centering
\includegraphics[width=\textwidth]{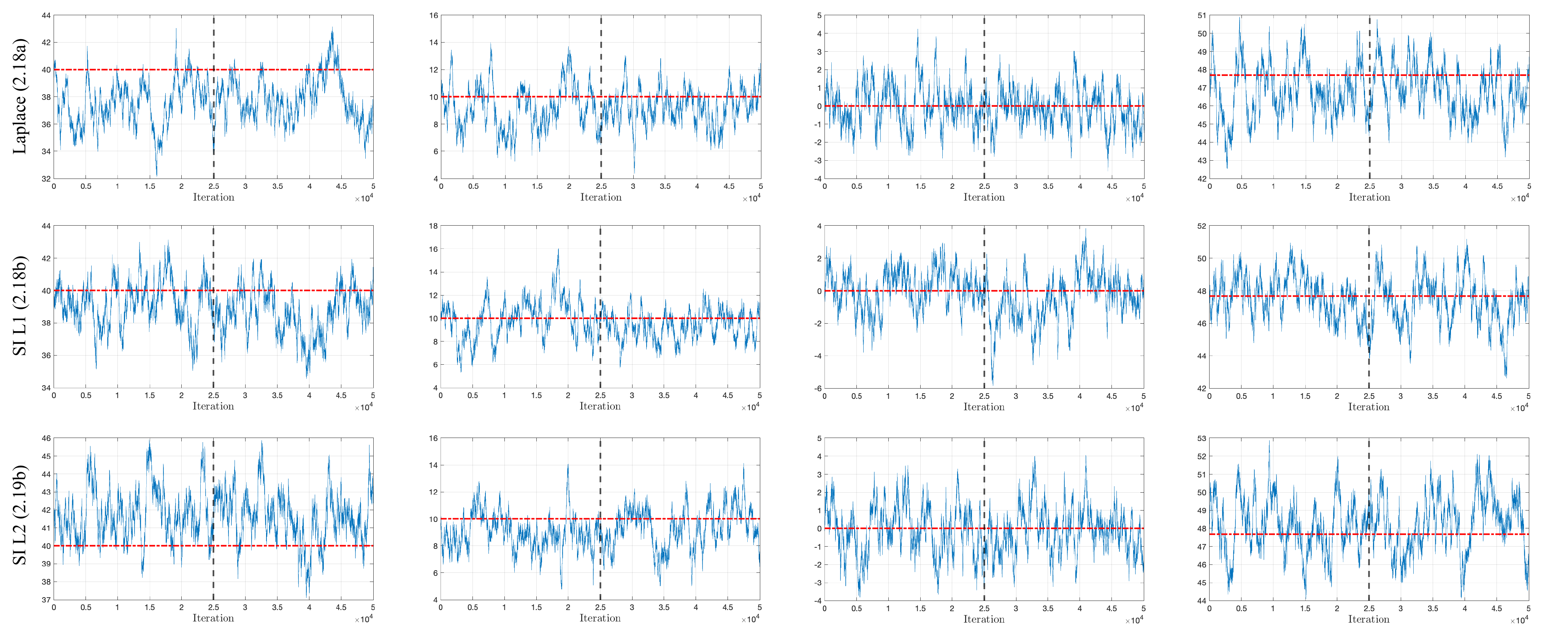}
\caption{Trace plots calculated at the points described by the red crosses in Figure \ref{fig:mask_demo}. Samples are taken from the posterior distributions defined using the (top)  Laplace prior (\ref{eq:PosteriorDensitylaplace}); (second row) support informed sparsity prior (\ref{eq:PosteriorDensitymask}); and (bottom) support informed $\ell_2$ prior (\ref{eq:PosteriorDensitymaskl2}).  Here $\sigma = 10$ ($SNR \approx -.2$ dB).  {The trace plots display similar behavior when $\sigma \approx 5.5$ ($SNR = 5$ dB).}} 
\label{fig:trace}
\end{figure}

Figure \ref{fig:trace} plots the value of each sample drawn from the posterior distribution {determined via Algorithm \ref{alg:MHAlgorithm}} described by (top) (\ref{eq:PosteriorDensitylaplace}), (middle) (\ref{eq:PosteriorDensitymask}), and (bottom) \eqref{eq:PosteriorDensitymaskl2} at the values in the domain indicated by the red crosses in Figure \ref{fig:mask_demo} {for $SNR = 5$ dB ($\sigma \approx 5.5$).}  {While each trace plot yields similar behavior, it is evident that at all points considered, using the posterior defined with the support informed sparsity prior (\ref{eq:PosteriorDensitymask}) results in samples that are more uniformly distributed about the true function value, indicated by the horizontal red line.}  {We note that comparable results were obtained for $\sigma = 10$ ($SNR \approx -.2$ dB) and are not displayed here.} 

\begin{figure}[htb]
\centering
\includegraphics[width=.5\textwidth]{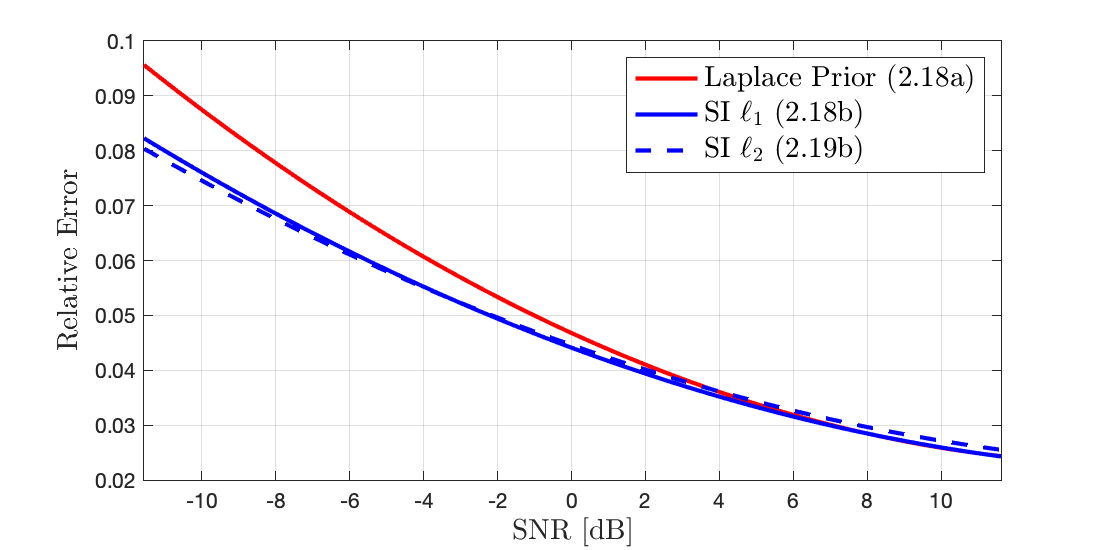}
\caption{The relative reconstruction error (\ref{eq:rel_error}) calculated across the entire spatial domain for decreasing SNR (\ref{eq:SNR}).} 
\label{fig:rel_error}
\end{figure}

{As was done for the sparse signal recovery case,} we also evaluate the accuracy of the reconstructions (\ref{eq:final_approx}) for Example \ref{ex:true_1D} given {discrete Fourier} measurements acquired with increasing levels of complex additive Gaussian white noise in (\ref{eq:ForwardModel}). Figure \ref{fig:rel_error} displays the relative error, \AG{\eqref{eq:rel_error}},
and Figure \ref{fig:abs_error} displays the absolute error, $|x_i-x^{true}_i|$, at the points in the spatial domain described by the red crosses in Figure \ref{fig:mask_demo}. {Each point $x_i$ across the domain is computed by averaging the results over 20 independent MCMC trials.}  
{Once again we see that} the support informed priors consistently yield either better or comparable results for all levels of SNR.  {This is true in smooth regions as well as near discontinuities, in which case it is clear that using the Laplace prior in regions where the underlying signal is {\em not} sparse in the presumably sparse domain yields dramatically worse results.}

%% file: sparsesignal.tex
{We first consider the recovery of a sparse signal. In this case the sparsifying transform matrix operator ${\mathcal L} = {\mathcal I}$ in \eqref{eq:final_posterior} and \eqref{eq:final_posteriorl2}.   This example eliminates any complications in evaluating the posterior  caused by different choices of sparsifying transform operators.   Moreover, sparse signal recovery is perhaps the simplest case to analyze, since conjugate priors are easily obtained, and indeed there is a closed form posterior for \eqref{eq:PosteriorDensityhierarchical}.\footnote{{As noted previously, the focus of our investigation is not to compare the efficiency of sampling methods for closed form posteriors. Hence we will still employ the MH algorithm even though Gibbs sampling is more appropriate in this particular case.}}  As noted previously, in the case of sparse signals, SBL with EM  yields a comparable posterior, {in the sense that it uses an empirically based  support informed prior}, to that in \eqref{eq:PosteriorDensitymaskl2}. It is not computationally efficient in higher dimensions, however, \cite{churchill2019detecting,tipping2001sparse,wipf2004sparse}.  To further simplify the analysis we choose the non-zero values of the underlying signal to be exactly $1$.}


\begin{figure}[htb]
\centering
\includegraphics[width = .8\textwidth]{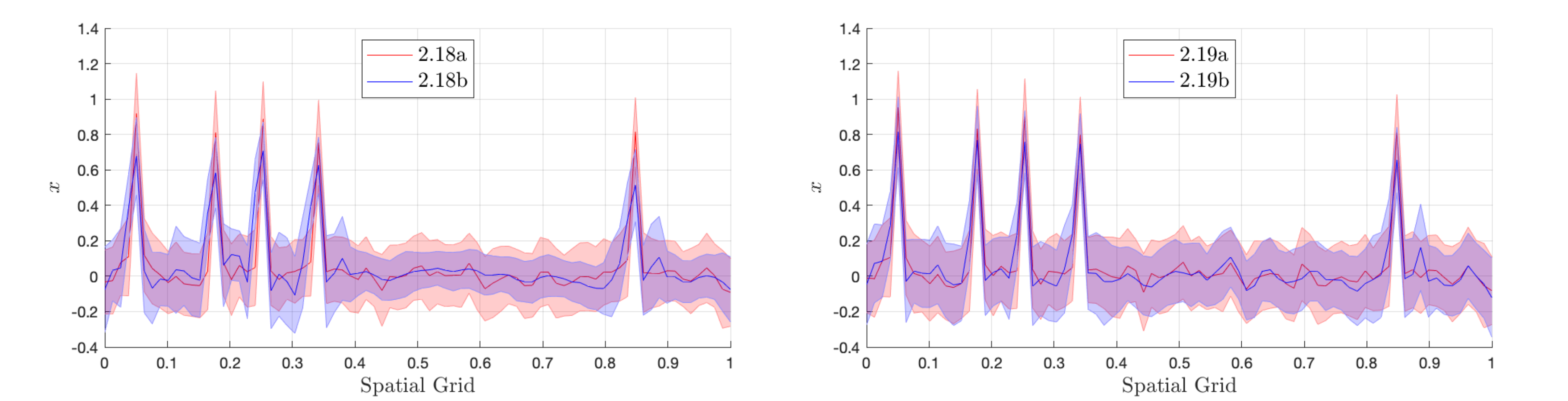}
\caption{{The mean and $C_{0.05}$ credibility intervals when the input data are noisy Fourier samples with $\text{SNR} = 0$dB. (left) comparison of \eqref{eq:PosteriorDensitylaplace} and  \eqref{eq:PosteriorDensitymask}; (right) comparison of \eqref{eq:PosteriorDensityhierarchical} and  \eqref{eq:PosteriorDensitymaskl2}.}}
\label{fig:sparse_sig_1D}
\end{figure} 
{Figure \ref{fig:sparse_sig_1D} displays the results and the $C_{0.05}$ credibility intervals where the observations in  \eqref{eq:ForwardModel} are noisy Fourier samples with $\text{SNR} = 0$dB, so that $A$ is modeled by the discrete Fourier transform,  and the complex additive white Gaussian noise is given by $\mathbf{E}\sim\mathcal{CN}(0,\sigma^2\mathbb{I}_{{n}})$ with $\sigma \approx 4.5$.      Other credibility intervals (e.g. $C_{0.01}$ and $C_{0.1}$) for $\text{SNR} = 0$dB produce comparable results, with reduced uncertainty especially apparent in the sparse regions of the signal when our new support informed prior is used.  {As shown respectively in Figure \ref{fig:sparse_sig_1D}(left) and (right), this is true for both \eqref{eq:PosteriorDensitymask} and \eqref{eq:PosteriorDensitymaskl2}, although the results are more striking for the sparse prior ($\ell_1$) case.}
\begin{figure}[htb]
\centering
\includegraphics[width = .8\textwidth]{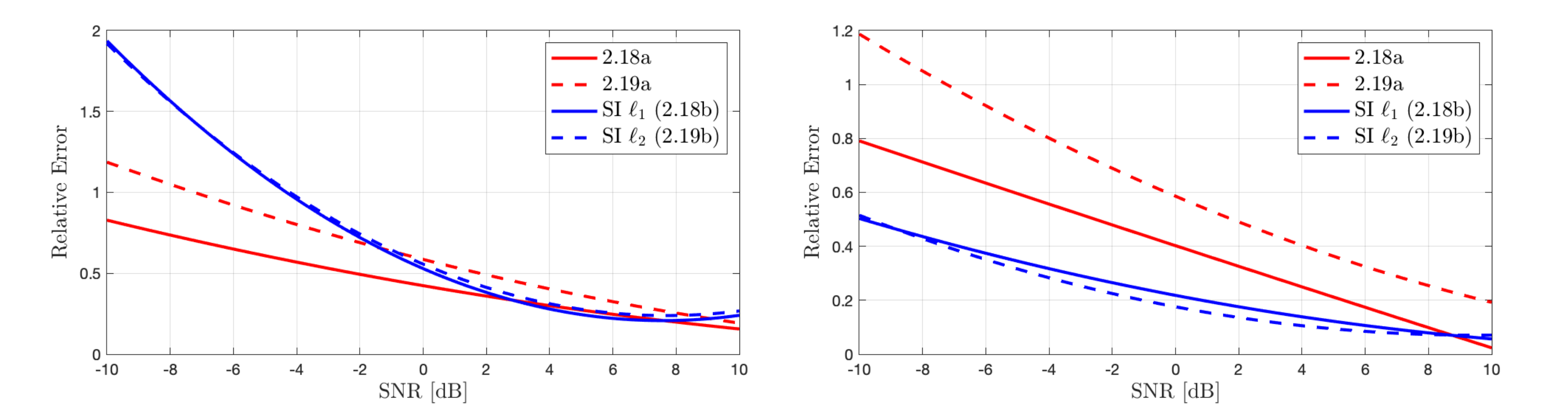}
\caption{{Relative error of the mean (\ref{eq:rel_error}) for $\text{SNR} \in [-10,10]$ at (left) a location where the sparse signal equals zero and (right) at a location where the sparse signal is non-zero.}} 
\label{fig:sparse_sig_1D_error}
\end{figure} 

{Figure \ref{fig:sparse_sig_1D_error}  compares the relative error of the mean (\ref{eq:rel_error}) using \eqref{eq:PosteriorDensitylaplace},  \eqref{eq:PosteriorDensitymask}, \eqref{eq:PosteriorDensityhierarchical} and  \eqref{eq:PosteriorDensitymaskl2} for decreasing $\text{SNR}$.  These comparisons are made at locations in the domain where the signal value is zero (left) and where it is $1$ (right).  From Figure \ref{fig:sparse_sig_1D_error}(right) we observe that by using the support informed prior, {that is by ``masking out''} the sparse prior at the non-zero function locations, the target density is equivalent to likelihood in regions of support, i.e.~we are not enforcing an {\em incorrect} sparse prior in those locations.  This is true for all levels of noise.  The relative error in sparse regions, depicted in Figure \ref{fig:sparse_sig_1D_error}{(left)}, demonstrates comparable results for $\text{SNR} \ge 0$.  Once the noise dominates the signal, the mean of the posterior does not yield meaningful results using either prior.}\footnote{{In this case it is probably better to employ \eqref{eq:PosteriorDensitylaplace} using a Gibbs sampler,} \cite{MCMCR}.} {For low $\text{SNR}$, the results using the standard (not support informed) priors are essentially the same regardless of whether or not the signal is zero or non-zero.  When combining these results with those in Figure \ref{fig:sparse_sig_1D}, we see that even though the support informed sparse prior may not be as accurate in the mean when the signal value is zero, the uncertainty is still reduced in those regions.  Thus we conclude that using the support informed sparse prior yields improved accuracy in the mean at the non-zero signal locations while also reducing the uncertainty everywhere.}

%% file: conclusions.tex
This paper developed a new empirical Bayesian inference method using joint sparsity. In particular, the joint sparsity of  multiple measurement vectors (MMVs) of observable data in the sparse domain is used to determine the likely support in the sparse domain.  From this information a new {\em support informed sparse prior} is developed, which is then used to generate a posterior and subsequently sampled using the MCMC algorithm.  The support in the sparse domain is obtained using provable convergence rates of edge detection methods from the observable data (in this case Fourier data) near and away from that support, \cite{archibaldgelb2002,archibald2005polynomial,GT1,GT06,VBJS2}.  

Our numerical experiments for sparse signals and piecewise smooth functions demonstrate that the new support informed sparse prior yields reduced uncertainty when compared to the Laplace prior, and improved accuracy in the mean wherever the underlying signal was non-zero in the sparse domain.  These findings 
strongly suggest that employing a {\em any} global (unmasked) sparse prior, regardless of the type, leads to less accurate results.  We also observed faster convergence in the ACF plots in smooth regions and near discontinuities.  This was more apparent for signals with higher SNR, but still true when $SNR \approx 0$ dB.

The support informed sparse prior does not yield a closed form solution for the posterior, so we therefore use the MH algorithm to sample it. This is of course not as efficient as Gibbs sampling. However, because {our method samples} at each individual point  in the domain, {meaning that separate chains are independently constructed for each pixel,} extensions to multi-dimensional images is straightforward. {Our approach may also be expanded  to include complex valued signals and images.  In this case, it was shown in \cite{ArchibaldGelb,churchill2018edge,sequentialjointsparsity,scarnati2018joint} how the PA transform may be efficiently constructed in multi-dimensions to promote sparsity in the real part of the respective signal or image.}   Finally, a main advantage of our technique is that once the support in the sparse domain is determined, different priors can be used in each smooth region, and in some cases a more efficient MAP estimate will suffice.  {As a main motivation for this investigation is to reduce the uncertainty in applications such as SAR imaging, we will explore these ideas} in future work.

%% file: appendixPA.tex
While a variety of sparsifying operators can be used for \eqref{eq:final_posterior}, for this investigation we employ the polynomial annihilation (PA) operator, ${\mathcal L} = {\mathcal L}^m$, which is based on the PA  edge detection method, \cite{archibald2005polynomial}.  The primary advantage in using the PA operator is that it is constructed to yield higher order approximations whenever $m > 1$. (When $m = 1$, it is equivalent to TV.)  This means that the approximation of the corresponding  edge function of a piecewise smooth function  is more likely to be sparse, and in particular yields values $\approx 0$ whenever the corresponding grid point  falls in a smooth region of the function.  This is especially important when the smooth region has a lot of variation (as opposed to being piecewise constant) or if the measurements are  sparsely sampled.  Below we describe the construction of  ${\mathcal L}^m$.  More information can be found in \cite{ArchibaldGelb}.

For ease of presentation and for consistency with our manuscript we consider $f(x)$ to be a piecewise smooth function on $[-\pi,\pi]$ (the finite domain is arbitrary) given on uniform grid point values $x_j = -\pi + \frac{2\pi j}{n}$, $j = 1,\cdots,n$.  We note that this is not a requirement for the PA method,  and furthermore that the method is designed to approximate the edge domain for $d > 1$ dimensional  piecewise smooth functions for any given set of scattered points of data, \cite{archibald2005polynomial}. However, since images are usually defined on a Cartesian grid,  for  the purposes of defining a sparsifying transform operator to be used in the prior given by (\ref{eq:laplace_prior}), as well as all of the priors defined in sequel throughout the manuscript, it was demonstrated in \cite{ArchibaldGelb} that applying the PA transform operator dimension by dimension was accurate and efficient.  

With this in mind we now define the $m$th order PA edge detection approximation as
\begin{equation}
\label{eq:EdgeDetector}
L^mf(x)=\frac{1}{q^{m}(x)}\sum_{x_j\in S_x}c_j(x)f(x_j),
\end{equation}
where $S_x$ is the local set of $m+1$ grid points from the set of given grid points about $x$, $c_j(x)$ are the polynomial annihilation edge detection coefficients, (\ref{eq:EdgeCoeffs}), and $q^m(x)$ is the normalization factor, (\ref{eq:normalization}).  Each parameter of the method can be further described as:
\begin{itemize}
\item $S_x$: For any particular cell $I_j=[x_j,x_{j+1})$, there are $m$ possible stencils, $S_x$ of size $m+1$, that contain the interval $I_j$.  For simplicity, we assume that the stencils are centered around the interval of interest, $I_j$, and are given by
$$S_{I_j} = \{x_{j - \frac{m}{2}},\cdots,x_{j+\frac{m}{2}}\},\hspace{.2in}
S_{I_j} = \{x_{j - \frac{m+1}{2}},\cdots,x_{j+\frac{m-1}{2}}\}$$
for $m$ even and odd respectively.  For non-periodic or non-zero-padded solutions the stencils are adapted to be more one sided as the boundaries of the interval are approached, \cite{archibald2005polynomial}.  To avoid cumbersome notation, we write $S_x$ as the generic stencil unless further clarification is needed.
\item $c_j(x)$:  The polynomial annihilation edge detection  coefficients, $c_j(x), \; j=1,\dots, m+1$, are constructed to annihilate polynomials up to degree $m$. They are obtained by solving the system
\begin{equation}
\label{eq:EdgeCoeffs}
\sum_{x_j\in S_x}c_j(x)p_\ell(x_j)= p_\ell^{(m)}(x), \; j=1,\dots, m+1,
\end{equation}
where $p_\ell$, $\ell=0,\dots,m,$ is a basis of for the space of polynomials of degree $\le m$.
\item $q^{m}(x)$:  The normalization factor, $q^m(x)$, normalizes the approximation to assure the proper convergence of $L^mf$ to the jump value at each discontinuity.  It is computed as
\begin{equation}
\label{eq:normalization}
q^m(x) = \sum_{x_j\in S_x^+} c_j(x),
\end{equation}
where $S_x^+$ is the set of points $x_j \in S_x$ such that $x_j \ge x$.
\end{itemize}

Since the data are given on uniform points,  there is an explicit formula for the PA edge detection coefficients, independent of location $x$, computed as (\cite{archibald2005polynomial})
\begin{equation}
\label{eq:UniformCoefficients}
c_j=\frac{m!}{\prod_{k=1, k \neq j}^{m+1}(j-k)\Delta x},\quad j=1,\dots ,m+1, \quad \Delta x = \frac{2\pi}{n}.
\end{equation}
We can now define the PA transform matrix ${\mathcal L}^m$ {used in the sparse prior} in terms of its elements as
\begin{equation}
\label{eq:PAtransform}
{\mathcal L}^m_{j,l}  = \frac{c(j,l)}{q^m(x_l)}, \hspace{.2in} 0 \le l \le n, \hspace{.2in} 0 \le j < n,
\end{equation}
where
$$c(j,l) = \left\{
     \begin{array}{lr}
       c_{j-l-\lfloor \frac{m}{2}\rfloor}, & 0<j-l-\lfloor \frac{m}{2}\rfloor+s(j,l)\leq m+1\\
       0 & \textrm{otherwise}
     \end{array}
   \right.
$$
and
$$s(j,l) = \left\{
     \begin{array}{lr}
       l-\lfloor \frac{m}{2}\rfloor, & \textrm{$l\leq\lfloor \frac{m}{2}\rfloor$}\\
        l+m-\lfloor \frac{m}{2}\rfloor-n, & \textrm{$l+m-\lfloor \frac{m}{2}\rfloor>n$}\\
       0 & \textrm{otherwise.}
     \end{array}
   \right.
$$
If the underlying image is known to be periodic, or is zero padded at the boundaries, a centered stencil can be used throughout the domain producing a circulant matrix ${\mathcal L}^m$.  A reduction of accuracy is expected near the boundaries for the non-periodic case due to the one-sided stencils.  

For example, assuming periodicity, the banded matrix ${\mathcal L}^1$  and ${\mathcal L}^3$ are given by
       \begin{equation}
        \label{eq:PAexamples}
                \mathcal{L}^1 = \begin{bmatrix}
                  1 & -1 &        &         &      \\
                    & 1  & -1     &         &      \\
                    &    & \ddots & \ddots  &       \\
                    &    &        &  1      & -1    \\
                 -1 &    &        &         & 1
                 \end{bmatrix}, \quad
                \mathcal{L}^3 = \frac{1}{2}\begin{bmatrix}
                  3  & -3     & 1      &        &                 & -1 \\
                  -1 &  3     & -3     & 1      &                 &     \\
                     & \ddots & \ddots & \ddots & \ddots          &     \\
                     &        &  1      &    3    &  -3           & 1   \\
                  1  &        &        & 1        &  3            & -3 \\
                 -3  & 1      &       &           &   -1          & 3
                 \end{bmatrix}.
        \end{equation}




%% file: Bayes_JS_SIAMUQ.bbl
\begin{thebibliography}{10}

\bibitem{VBJSGelb}
{\sc B.~Adcock, A.~Gelb, G.~Song, and Y.~Sui}, {\em Joint sparse recovery based
  on variances}, SIAM Journal on Scientific Computing, 41 (2019),
  pp.~A246--A268, \url{https://doi.org/10.1137/17M1155983},
  \url{https://doi.org/10.1137/17M1155983},
  \url{https://arxiv.org/abs/https://doi.org/10.1137/17M1155983}.

\bibitem{archibaldgelb2002}
{\sc R.~Archibald and A.~Gelb}, {\em Reducing the effects of noise in image
  reconstruction}, Journal of Scientific Computing, 17 (2002), pp.~167--180.

\bibitem{ArchibaldGelb}
{\sc R.~Archibald, A.~Gelb, and R.~Platte}, {\em Image reconstruction from
  undersampled fourier data using the polynomial annihilation transform},
  Journal of Scientific Computing,  (2015),
  \url{https://doi.org/10.1007/s10915-015-0088-2}.

\bibitem{archibald2005polynomial}
{\sc R.~Archibald, A.~Gelb, and J.~Yoon}, {\em Polynomial fitting for edge
  detection in irregularly sampled signals and images}, SIAM Journal on
  Numerical Analysis, 43 (2005), pp.~259--279.

\bibitem{BardsleyUQ}
{\sc J.~M. Bardsley}, {\em Mcmc-based image reconstruction with uncertainty
  quantification}, SIAM Journal on Scientific Computing, 34 (2012),
  pp.~A1316--A1332, \url{https://doi.org/10.1137/11085760X},
  \url{https://doi.org/10.1137/11085760X},
  \url{https://arxiv.org/abs/https://doi.org/10.1137/11085760X}.

\bibitem{bhadra2019lasso}
{\sc A.~Bhadra, J.~Datta, N.~G. Polson, B.~Willard, et~al.}, {\em Lasso meets
  horseshoe: A survey}, Statistical Science, 34 (2019), pp.~405--427.

\bibitem{calvetti2020hybrid}
{\sc D.~Calvetti, M.~Pragliola, and E.~Somersalo}, {\em Hybrid solver for
  hierarchical {B}ayesian inverse problems}, arXiv preprint arXiv:2003.06532,
  (2020).

\bibitem{calvetti2020sparse}
{\sc D.~Calvetti, M.~Pragliola, E.~Somersalo, and A.~Strang}, {\em Sparse
  reconstructions from few noisy data: analysis of hierarchical {B}ayesian
  models with generalized gamma hyperpriors}, Inverse Problems, 36 (2020),
  p.~025010.

\bibitem{Candes}
{\sc E.~J. {Candes}, J.~{Romberg}, and T.~{Tao}}, {\em Robust uncertainty
  principles: exact signal reconstruction from highly incomplete frequency
  information}, IEEE Transactions on Information Theory, 52 (2006),
  pp.~489--509.

\bibitem{candes2008enhancing}
{\sc E.~J. Candes, M.~B. Wakin, and S.~P. Boyd}, {\em Enhancing sparsity by
  reweighted $\ell_1$ minimization}, Journal of {F}ourier analysis and
  applications, 14 (2008), pp.~877--905.

\bibitem{CompressiveSensingCandes}
{\sc E.~J. Candès, J.~K. Romberg, and T.~Tao}, {\em Stable signal recovery
  from incomplete and inaccurate measurements}, Communications on Pure and
  Applied Mathematics, 59 (2006), pp.~1207--1223,
  \url{https://doi.org/10.1002/cpa.20124},
  \url{https://onlinelibrary.wiley.com/doi/abs/10.1002/cpa.20124},
  \url{https://arxiv.org/abs/https://onlinelibrary.wiley.com/doi/pdf/10.1002/cpa.20124}.

\bibitem{casella2002statistical}
{\sc G.~Casella and R.~Berger}, {\em Statistical Inference}, Duxbury advanced
  series, Duxbury Thomson Learning, 2002,
  \url{https://books.google.com/books?id=ZpkPPwAACAAJ}.

\bibitem{MMVChenHuo}
{\sc J.~{Chen} and X.~{Huo}}, {\em Theoretical results on sparse
  representations of multiple-measurement vectors}, IEEE Transactions on Signal
  Processing, 54 (2006), pp.~4634--4643.

\bibitem{churchill2018edge}
{\sc V.~Churchill, R.~Archibald, and A.~Gelb}, {\em Edge-adaptive $\ell_2$
  regularization image reconstruction from non-uniform {F}ourier data}, arXiv
  preprint arXiv:1811.08487,  (2018).
\newblock https://arxiv.org/abs/1811.08487.

\bibitem{churchill2019detecting}
{\sc V.~Churchill and A.~Gelb}, {\em Detecting edges from non-uniform {F}ourier
  data via sparse {B}ayesian learning}, Journal of Scientific Computing,
  (2019), pp.~1--22.

\bibitem{daubechies2010iteratively}
{\sc I.~Daubechies, R.~DeVore, M.~Fornasier, and C.~S. G{\"u}nt{\"u}rk}, {\em
  Iteratively reweighted least squares minimization for sparse recovery},
  Communications on Pure and Applied Mathematics, 63 (2010), pp.~1--38.

\bibitem{dempster1977maximum}
{\sc A.~P. Dempster, N.~M. Laird, and D.~B. Rubin}, {\em Maximum likelihood
  from incomplete data via the em algorithm}, Journal of the royal statistical
  society. Series B (methodological),  (1977), pp.~1--38.

\bibitem{Donoho}
{\sc D.~L. {Donoho}}, {\em Compressed sensing}, IEEE Transactions on
  Information Theory, 52 (2006), pp.~1289--1306.

\bibitem{ConvergenceMCMC}
{\sc J.~Ellis}, {\em A practical guide to mcmc part 1: Mcmc basics}, 2018,
  \url{https://jellis18.github.io/post/2018-01-02-mcmc-part1/} (accessed
  2020-04-29).

\bibitem{MCMC}
{\sc D.~Gamerman and H.~F. Lopes}, {\em Morkov chain monte carlo}, 68 (2006).
\newblock MCMC.

\bibitem{VBJS}
{\sc A.~Gelb and T.~Scarnati}, {\em Reducing effects of bad data using variance
  based joint sparsity recovery}, Journal of Scientific Computing, 78 (2019),
  pp.~94--120.

\bibitem{GT1}
{\sc A.~Gelb and E.~Tadmor}, {\em Detection of edges in spectral data}, Applied
  and Computational Harmonic Analysis, 7 (1999), pp.~101--135.

\bibitem{GT06}
{\sc A.~Gelb and E.~Tadmor}, {\em Adaptive edge detectors for piecewise smooth
  data based on the minmod limiter}, J. Sci. Comput., 28 (2006), pp.~279--306.

\bibitem{hastings}
{\sc W.~K. Hastings}, {\em {Monte Carlo sampling methods using Markov chains
  and their applications}}, Biometrika, 57 (1970), pp.~97--109,
  \url{https://doi.org/10.1093/biomet/57.1.97},
  \url{https://doi.org/10.1093/biomet/57.1.97},
  \url{https://arxiv.org/abs/https://academic.oup.com/biomet/article-pdf/57/1/97/23940249/57-1-97.pdf}.

\bibitem{kaipio2006statistical}
{\sc J.~Kaipio and E.~Somersalo}, {\em Statistical and Computational Inverse
  Problems}, Applied Mathematical Sciences, Springer New York, 2006,
  \url{https://books.google.com/books?id=h0i-Gi4rCZIC}.

\bibitem{metropolis}
{\sc N.~Metropolis, A.~W. Rosenbluth, M.~N. Rosenbluth, A.~H. Teller, and
  E.~Teller}, {\em Equation of state calculations by fast computing machines},
  The Journal of Chemical Physics, 21 (1953), pp.~1087--1092,
  \url{https://doi.org/10.1063/1.1699114},
  \url{https://doi.org/10.1063/1.1699114},
  \url{https://arxiv.org/abs/https://doi.org/10.1063/1.1699114}.

\bibitem{murphy2012machine}
{\sc K.~Murphy}, {\em Machine Learning: A Probabilistic Perspective}, Adaptive
  Computation and Machine Learning series, MIT Press, 2012,
  \url{https://books.google.com/books?id=NZP6AQAAQBAJ}.

\bibitem{piironen2017sparsity}
{\sc J.~Piironen, A.~Vehtari, et~al.}, {\em Sparsity information and
  regularization in the horseshoe and other shrinkage priors}, Electronic
  Journal of Statistics, 11 (2017), pp.~5018--5051.

\bibitem{acf_matlab}
{\sc C.~Price}, {\em Autocorrelation function {(ACF)}}.
\newblock
  \url{https://www.mathworks.com/matlabcentral/fileexchange/30540-autocorrelation-function-acf},
  2020.
\newblock accessed: October 2020.

\bibitem{MCMCR}
{\sc C.~Robert and G.~Casella}, {\em Introducing Monte Carlo Methods with R},
  Use R!, Springer, 2010, \url{https://books.google.com/books?id=WIjMyiEiHCsC}.

\bibitem{sequentialjointsparsity}
{\sc T.~Sanders, A.~Gelb, and R.~Platte}, {\em Composite sar imaging using
  sequential joint sparsity}, Journal of Computational Physics, 338 (2017),
  \url{https://doi.org/10.1016/j.jcp.2017.02.071}.

\bibitem{scarnati2018joint}
{\sc T.~Scarnati and A.~Gelb}, {\em Joint image formation and two-dimensional
  autofocusing for synthetic aperture radar data}, Journal of Computational
  Physics, 374 (2018), pp.~803--821.

\bibitem{VBJS2}
{\sc T.~Scarnati and A.~Gelb}, {\em Accurate and efficient image reconstruction
  from multiple measurements of fourier samples}, Journal of Computational
  Mathematics, 88 (2020), pp.~798--828.

\bibitem{tipping2001sparse}
{\sc M.~E. Tipping}, {\em Sparse bayesian learning and the relevance vector
  machine}, Journal of machine learning research, 1 (2001), pp.~211--244.

\bibitem{wipf2004sparse}
{\sc D.~P. Wipf and B.~D. Rao}, {\em Sparse bayesian learning for basis
  selection}, IEEE Transactions on Signal processing, 52 (2004),
  pp.~2153--2164.

\end{thebibliography}
